\catcode`\@=11

\font\tenmsa=msam10 \font\sevenmsa=msam7 \font\fivemsa=msam5
\font\tenmsb=msbm10
\font\sevenmsb=msbm7 \font\fivemsb=msbm5 \newfam\msafam
\newfam\msbfam
\textfont\msafam=\tenmsa \scriptfont\msafam=\sevenmsa
\scriptscriptfont\msafam=\fivemsa \textfont\msbfam=\tenmsb
\scriptfont\msbfam=\sevenmsb \scriptscriptfont\msbfam=\fivemsb

\def\hexnumber@#1{\ifnum#1<10 \number#1\else \ifnum#1=10
A\else\ifnum#1=11
 B\else\ifnum#1=12 C\else \ifnum#1=13 D\else\ifnum#1=14
E\else\ifnum#1=15
 F\fi\fi\fi\fi\fi\fi\fi}

\def\msa@{\hexnumber@\msafam} \def\msb@{\hexnumber@\msbfam}
\mathchardef\boxdot="2\msa@00 \mathchardef\boxplus="2\msa@01
\mathchardef\boxtimes="2\msa@02 \mathchardef\square="0\msa@03
\mathchardef\blacksquare="0\msa@04 \mathchardef\centerdot="2\msa@05
\mathchardef\lozenge="0\msa@06 \mathchardef\blacklozenge="0\msa@07
\mathchardef\circlearrowright="3\msa@08
\mathchardef\circlearrowleft="3\msa@09
\mathchardef\rightleftharpoons="3\msa@0A
\mathchardef\leftrightharpoons="3\msa@0B
\mathchardef\boxminus="2\msa@0C
\mathchardef\Vdash="3\msa@0D \mathchardef\Vvdash="3\msa@0E
\mathchardef\vDash="3\msa@0F \mathchardef\twoheadrightarrow="3\msa@10
\mathchardef\twoheadleftarrow="3\msa@11
\mathchardef\leftleftarrows="3\msa@12
\mathchardef\rightrightarrows="3\msa@13
\mathchardef\upuparrows="3\msa@14
\mathchardef\downdownarrows="3\msa@15
\mathchardef\upharpoonright="3\msa@16

\mathchardef\downharpoonright="3\msa@17
\mathchardef\upharpoonleft="3\msa@18
\mathchardef\downharpoonleft="3\msa@19
\mathchardef\rightarrowtail="3\msa@1A
\mathchardef\leftarrowtail="3\msa@1B
\mathchardef\leftrightarrows="3\msa@1C
\mathchardef\rightleftarrows="3\msa@1D
\mathchardef\Lsh="3\msa@1E \mathchardef\Rsh="3\msa@1F
\mathchardef\rightsquigarrow="3\msa@20
\mathchardef\leftrightsquigarrow="3\msa@21
\mathchardef\looparrowleft="3\msa@22
\mathchardef\looparrowright="3\msa@23 \mathchardef\circeq="3\msa@24
\mathchardef\succsim="3\msa@25 \mathchardef\gtrsim="3\msa@26
\mathchardef\gtrapprox="3\msa@27 \mathchardef\multimap="3\msa@28
\mathchardef\therefore="3\msa@29 \mathchardef\because="3\msa@2A
\mathchardef\doteqdot="3\msa@2B 
\mathchardef\traceiangleq="3\msa@2C \mathchardef\precsim="3\msa@2D
\mathchardef\lesssim="3\msa@2E \mathchardef\lessapprox="3\msa@2F
\mathchardef\eqslantless="3\msa@30 \mathchardef\eqslantgtr="3\msa@31
\mathchardef\curlyeqprec="3\msa@32 \mathchardef\curlyeqsucc="3\msa@33
\mathchardef\preccurlyeq="3\msa@34 \mathchardef\leqq="3\msa@35
\mathchardef\leqslant="3\msa@36 \mathchardef\lessgtr="3\msa@37
\mathchardef\backprime="0\msa@38 \mathchardef\risingdotseq="3\msa@3A
\mathchardef\fallingdotseq="3\msa@3B
\mathchardef\succcurlyeq="3\msa@3C
\mathchardef\geqq="3\msa@3D \mathchardef\geqslant="3\msa@3E
\mathchardef\gtrless="3\msa@3F \mathchardef\sqsubset="3\msa@40
\mathchardef\sqsupset="3\msa@41
\mathchardef\trianglerighteq="3\msa@44
\mathchardef\trianglelefteq="3\msa@45 \mathchardef\bigstar="0\msa@46
\mathchardef\between="3\msa@47
\mathchardef\blacktriangledown="0\msa@48
\mathchardef\blacktriangleright="3\msa@49
\mathchardef\blacktriangleleft="3\msa@4A
\mathchardef\blacktriangle="0\msa@4E
\mathchardef\triangledown="0\msa@4F
\mathchardef\eqcirc="3\msa@50 \mathchardef\lesseqgtr="3\msa@51
\mathchardef\gtreqless="3\msa@52 \mathchardef\lesseqqgtr="3\msa@53
\mathchardef\gtreqqless="3\msa@54 \mathchardef\Rrightarrow="3\msa@56
\mathchardef\Lleftarrow="3\msa@57 \mathchardef\veebar="2\msa@59
\mathchardef\barwedge="2\msa@5A \mathchardef\doublebarwedge="2\msa@5B
\mathchardef\angle="0\msa@5C \mathchardef\measuredangle="0\msa@5D
\mathchardef\sphericalangle="0\msa@5E
\mathchardef\varpropto="3\msa@5F
\mathchardef\smallsmile="3\msa@60 \mathchardef\smallfrown="3\msa@61
\mathchardef\Subset="3\msa@62 \mathchardef\Supset="3\msa@63
\mathchardef\Cup="2\msa@64 
\mathchardef\Cap="2\msa@65
 \mathchardef\curlywedge="2\msa@66
\mathchardef\curlyvee="2\msa@67 \mathchardef\leftthreetimes="2\msa@68
\mathchardef\rightthreetimes="2\msa@69
\mathchardef\subseteqq="3\msa@6A
\mathchardef\supseteqq="3\msa@6B \mathchardef\bumpeq="3\msa@6C
\mathchardef\Bumpeq="3\msa@6D \mathchardef\lll="3\msa@6E

\mathchardef\ggg="3\msa@6F 
\mathchardef\circledS="0\msa@73
\mathchardef\pitchfork="3\msa@74 \mathchardef\dotplus="2\msa@75
\mathchardef\backsim="3\msa@76 \mathchardef\backsimeq="3\msa@77
\mathchardef\complement="0\msa@7B \mathchardef\intercal="2\msa@7C
\mathchardef\circledcirc="2\msa@7D \mathchardef\circledast="2\msa@7E
\mathchardef\circleddash="2\msa@7F
\def\ulcorner{\delimiter"4\msa@70\msa@70 }
\def\urcorner{\delimiter"5\msa@71\msa@71 }
\def\llcorner{\delimiter"4\msa@78\msa@78 }
\def\lrcorner{\delimiter"5\msa@79\msa@79 }
\def\yen{\mathhexbox\msa@55 }
\def\checkmark{\mathhexbox\msa@58 } \def\circledR{\mathhexbox\msa@72
}
\def\maltese{\mathhexbox\msa@7A } \mathchardef\lvertneqq="3\msb@00
\mathchardef\gvertneqq="3\msb@01 \mathchardef\nleq="3\msb@02
\mathchardef\ngeq="3\msb@03 \mathchardef\nless="3\msb@04
\mathchardef\ngtr="3\msb@05 \mathchardef\nprec="3\msb@06
\mathchardef\nsucc="3\msb@07 \mathchardef\lneqq="3\msb@08
\mathchardef\gneqq="3\msb@09 \mathchardef\nleqslant="3\msb@0A
\mathchardef\ngeqslant="3\msb@0B \mathchardef\lneq="3\msb@0C
\mathchardef\gneq="3\msb@0D \mathchardef\npreceq="3\msb@0E
\mathchardef\nsucceq="3\msb@0F \mathchardef\precnsim="3\msb@10
\mathchardef\succnsim="3\msb@11 \mathchardef\lnsim="3\msb@12
\mathchardef\gnsim="3\msb@13 \mathchardef\nleqq="3\msb@14
\mathchardef\ngeqq="3\msb@15 \mathchardef\precneqq="3\msb@16
\mathchardef\succneqq="3\msb@17 \mathchardef\precnapprox="3\msb@18
\mathchardef\succnapprox="3\msb@19 \mathchardef\lnapprox="3\msb@1A
\mathchardef\gnapprox="3\msb@1B \mathchardef\nsim="3\msb@1C
\mathchardef\napprox="3\msb@1D
\mathchardef\nsubseteqq="3\msb@22
\mathchardef\nsupseteqq="3\msb@23 \mathchardef\subsetneqq="3\msb@24
\mathchardef\supsetneqq="3\msb@25
\mathchardef\subsetneq="3\msb@28
\mathchardef\supsetneq="3\msb@29 \mathchardef\nsubseteq="3\msb@2A
\mathchardef\nsupseteq="3\msb@2B \mathchardef\nparallel="3\msb@2C
\mathchardef\nmid="3\msb@2D \mathchardef\nshortmid="3\msb@2E
\mathchardef\nshortparallel="3\msb@2F \mathchardef\nvdash="3\msb@30
\mathchardef\nVdash="3\msb@31 \mathchardef\nvDash="3\msb@32
\mathchardef\nVDash="3\msb@33 \mathchardef\ntrianglerighteq="3\msb@34
\mathchardef\ntrianglelefteq="3\msb@35
\mathchardef\ntriangleleft="3\msb@36
\mathchardef\ntriangleright="3\msb@37
\mathchardef\nleftarrow="3\msb@38
\mathchardef\nrightarrow="3\msb@39 \mathchardef\nLeftarrow="3\msb@3A
\mathchardef\nRightarrow="3\msb@3B
\mathchardef\nLeftrightarrow="3\msb@3C
\mathchardef\nleftrightarrow="3\msb@3D
\mathchardef\divideontimes="2\msb@3E
\mathchardef\varnothing="0\msb@3F \mathchardef\nexists="0\msb@40
\mathchardef\mho="0\msb@66 \mathchardef\thorn="0\msb@67
\mathchardef\beth="0\msb@69 \mathchardef\gimel="0\msb@6A
\mathchardef\daleth="0\msb@6B \mathchardef\lessdot="3\msb@6C
\mathchardef\gtrdot="3\msb@6D \mathchardef\ltimes="2\msb@6E
\mathchardef\rtimes="2\msb@6F \mathchardef\shortmid="3\msb@70
\mathchardef\shortparallel="3\msb@71
\mathchardef\smallsetminus="2\msb@72
\mathchardef\thicksim="3\msb@73 \mathchardef\thickapprox="3\msb@74
\mathchardef\approxeq="3\msb@75 \mathchardef\succapprox="3\msb@76
\mathchardef\precapprox="3\msb@77
\mathchardef\curvearrowleft="3\msb@78
\mathchardef\curvearrowright="3\msb@79 \mathchardef\digamma="0\msb@7A
\mathchardef\varkappa="0\msb@7B \mathchardef\hslash="0\msb@7D
\mathchardef\hbar="0\msb@7E \mathchardef\backepsilon="3\msb@7F
\def\Bbb{\ifmmode\let\next\Bbb@\else
\def\next{\errmessage{Use \string\Bbb\space only in math
mode}}\fi\next}
\def\Bbb@#1{{\Bbb@@{#1}}} \def\Bbb@@#1{\fam\msbfam#1}

\catcode`\@=\active

\def\CM{\hbox{{$\cal M$}}}

\def\cg{\hbox{{\sl g}}} 
\def\cb{\hbox{{\sl b}}} 
\def\cc{\hbox{{\sl c}}} 
\def\cf{\hbox{{\sl f}}}

\def\lform{\hbox{$\sqcup$}\llap{\hbox{$\sqcap$}}}

\def\h{{{1\over2}}}

\def\R{{\Bbb R}}
\def\C{{\Bbb C}}
\def\Z{{\Bbb Z}}

\def\eps{{\epsilon}}

\def\lcross{{>\!\!\!\triangleleft}}

\def\lcocross{{>\!\!\blacktriangleleft}}

\def\codcross{{\blacktriangleright\!\!\blacktriangleleft}}
\def\rbiprod{{\cdot\kern-.33em\triangleright\!\!\!<}}
\def\lbiprod{{>\!\!\!\triangleleft\kern-.33em\cdot\, }}

\def\tens{\mathop{\otimes}}

\def\la{{\triangleright}}\def\ra{{\triangleleft}}

\def\extd{{{\rm d}}}

\def\swap{{\leftrightarrow}}

\def\isom{{\cong}}

\def\Ad{{\rm Ad}}
\def\ad{{\rm ad}}

\def\ev{{\rm ev}}

\def\id{{\rm id}}

\def\<{\langle}
\def\>{\rangle}
\def\dila{{\varsigma}}

\def\equad{\kern -1.7em}

\def\eqn#1#2{\begin{equation}#2\label{#1}\end{equation}}

\def\o{{}_{\scriptscriptstyle(1)}}
\def\t{{}_{\scriptscriptstyle(2)}}

\def\bo{{}^{\bar{\scriptscriptstyle(1)}}}
\def\bt{{}^{\bar{\scriptscriptstyle(2)}}}

\def\und#1{{\underline {#1}}}

\def\uo{{{}^{\scriptscriptstyle(1)}}}
\def\ut{{{}^{\scriptscriptstyle(2)}}}

\def\Bo{{{}_{\und{\scriptscriptstyle(1)}}}}
\def\Bt{{{}_{\und{\scriptscriptstyle(2)}}}}

\def\text#1{\mbox{\rm #1}}
\def\note#1{}

\def\blacksquare{{\lform}}
\def\frac#1#2{{{#1\over#2}}}

\def\proof{\goodbreak\noindent{\bf Proof\quad}}

\def\endproof{{\ $\lform$}\bigskip }

\def\align#1{\begin{eqnarray*}#1\end{eqnarray*}}

\def\cmath#1{\[\begin{array}{c} #1 \end{array}\]}

\documentstyle[11pt]{article}
\newtheorem{lemma}{Lemma}[section]
\newtheorem{propos}[lemma]{Proposition}
\newtheorem{example}[lemma]{Example}
\newtheorem{theorem}[lemma]{Theorem}

\newtheorem{corol}[lemma]{Corollary}
\newtheorem{defin}[lemma]{Definition}

\begin{document}\baselineskip 20pt

{\ }\qquad\qquad \hskip 4.3in  DAMTP/96-98
\vspace{.2in}

\begin{center} {\LARGE BRAIDED LIE BIALGEBRAS}
\\ \baselineskip 13pt{\ }
{\ }\\
S. Majid\footnote{Royal Society University Research Fellow and Fellow
of
Pembroke College, Cambridge}\\
{\ }\\
Department of Mathematics, Harvard University\\
Science Center, Cambridge MA 02138, USA\footnote{During the calendar
years 1995 + 1996}\\
+\\
Department of Applied Mathematics \& Theoretical Physics\\
University of Cambridge, Cambridge CB3 9EW\\
\end{center}
\begin{center}
December 1996 -- revised February 1997
\end{center}

\vspace{20pt}
\begin{quote}\baselineskip 13pt
\noindent{\bf Abstract} We introduce braided Lie bialgebras as the
infinitesimal version of braided groups. They are Lie algebras and Lie
coalgebras with the coboundary of the Lie cobracket an infinitesimal
braiding. We provide theorems of transmutation, Lie biproduct, bosonisation and
double-bosonisation relating braided Lie bialgebras to usual Lie bialgebras.
Among the results, the kernel of any split projection of Lie bialgebras is a
braided-Lie bialgebra. The  Kirillov-Kostant Lie cobracket provides a natural
braided-Lie bialgebra on any complex simple Lie algebra $g$, as the
transmutation of the Drinfeld-Sklyanin Lie cobracket. Other
nontrivial braided-Lie bialgebras are associated to the inductive construction
of simple Lie bialgebras
along the $C$ and exceptional series.

\end{quote}
\baselineskip 23pt

\section{Introduction}

Braided geometry has been developed in recent years as a natural generalisation
of super-geometry with the role of $\Z/2\Z$ grading played by braid statistics.
It is also the kind of noncommutative geometry appropriate to quantum group
symmetry because the modules over a strict quantum group (a quasitriangular
Hopf algebra\cite{Dri}) form a braided category, hence any object covariant
under the quantum group is naturally braided. In particular, one has braided
groups\cite{Ma:bg} as generalisations of super-groups or super-Hopf algebras.
The famous quantum-braided plane with relations $yx=qxy$ is a braided group
with additive coproduct\cite{Ma:poi}.  We refer to
\cite{Ma:introm}\cite{Ma:book} for introductions to the 50-60 papers in which
the theory of braided groups is developed.

In a different direction, Drinfeld\cite{Dri:ham} has introduced Lie bialgebras
as an infinitesimalisation of the theory quantum groups. This concept has led
(on exponentiation) to an extensive theory of Poisson-Lie groups, as well as to
a Yang-Baxter theoretic approach to classical results of Lie theory, such as a
new proof of the Iwasawa decomposition; see \cite{Dri}\cite{Ma:book} for
reviews. We now combine these ideas by introducing the infinitesimal theory of
braided groups.  All computations and results will be in the setting of Lie
algebras, although motivated from the theory of braided groups.

In fact, there are several different concepts of precisely what one may mean by
the infinitesimal theory of braided groups. Firstly, one may keep the braided
category in which one works fixed and look at algebras which depart
infinitesimally from being commutative. In the category of vector spaces this
leads to Drinfeld's notion of Poisson-Lie group. Then one can consider the
coalgebra also in an infinitesimal form, which leads in the category of vector
spaces to Drinfeld's notion of Lie bialgebra. In the case of a braided category
one already has the notion of braided-Lie algebra\cite{Ma:lie} and, adding to
this,  one could similarly consider a Lie bialgebra in a braided category.
By contrast, we now go further and let the braiding also depart infinitesimally
from the usual vector space transposition. In principle, the degree of braiding
is independent of the degree of algebra commutativity or coalgebra
cocommutativity. Thus one could have infinitesimally braided algebras,
coalgebras and Hopf algebras as well. However, the case which appears to be of
most interest, on which we concentrate, is the case in which all three aspects
are made infinitesimal simultaneously, which we call a braided-Lie bialgebra.
The formal definition appears in Section~2. It consists of a Lie algebra $\cb$
equipped with further structure.

In Section~3 we provide the Lie version of the basic theorems from the theory
of braided groups. These basic theorems connect braided groups and quantum
groups by transmutation\cite{Ma:tra}\cite{Ma:bg} and
bosonisation\cite{Ma:bos}\cite{Ma:skl} procedures, thereby establishing (for
example) the existence of braided groups associated to all simple Lie algebras.
The theorems in Section~3 likewise connect braided-Lie bialgebras with
quasitriangular Lie bialgebras and establish the existence of the former. The
Lie versions of biproducts\cite{Rad:str} and of the more recent
double-bosonisation theorem \cite{Ma:dbos} are covered as well. For example,
the Lie version of the theory of biproducts states that the kernel of any split
Lie bialgebra projection $\cg\to\cf$ is a braided-Lie bialgebra $\cb$, and
$\cg=\cb\lbiprod\cf$.

In Section~4 we study some concrete examples of braided-Lie bialgebras,
including ones not obtained by transmutation. The simplest are ones with zero
braided-Lie cobracket as the infinitesimal versions of the q-affine plane
braided groups in \cite{Ma:poi}. As an application of braided-Lie bialgebras,
their bosonisations provide maximal parabolic or inhomogeneous Lie bialgebras.
Meanwhile, double-bosonisation allows the formulation in a basis-free way of
the notion of adjoining a node to a Dynkin diagram. For every simple Lie
bialgebra $\cg$ and braided-Lie bialgebra $\cb$ in its category of modules we
obtain a new simple Lie bialgebra $\cb\lbiprod\cg\rbiprod\cb^{*\rm op}$ as its
double-bosonisation This provides the inductive construction of all complex
simple Lie algebras,  complete with their Drinfeld-Sklyanin quasitriangular Lie
bialgebra structure (which is built up inductively at the same time). Some
concrete examples are given in detail.

These results have been briefly announced in \cite[Sec. 3]{Ma:asc}, of which
the present paper is the extended text. We work over a general ground field $k$
of characteristic not 2.

\subsection*{Acknowledgements} The results were obtained during a visit in June
1996 to the Mathematics Dept., Basel, Suisse. I thank the chairman for access
to the facilities.

\section{Braided-Lie bialgebras}

We will be concerned throughout with the Lie version of braided categories
obtained as module categories over quantum groups. In principle one could also
formulate an abstract notion of `infinitesimal braiding' as a Lie version of a
general braided category, but since no examples other than the ones related to
quantum groups are known we limit ourselves essentially to this concrete
setting. Some slight extensions (such as to Lie crossed modules) will be
considered as well, later on.

As a Lie version of a strict quantum group we use Drinfeld's notion of a
quasitriangular Lie bialgebra\cite{Dri:ham}\cite{Dri}. We recall that a Lie
bialgebra is a Lie algebra $\cg$ equipped with linear map $\delta:\cg\to
\cg\tens\cg$ forming a Lie coalgebra (in the finite dimensional case this is
equivalent to a Lie bracket on $\cg^*$) and being a 1-cocycle with values in
$\cg\tens\cg$ as a $\cg$-module by the natural extension of $\ad$. It is
quasitriangular if there exists $r\in\cg\tens\cg$ obeying $\extd r=\delta$ in
the Lie algebra complex, obeying the Classical Yang-Baxter Equation (CYBE)
\eqn{CYBE}{{}[r\uo,r'\uo]\tens r\ut\tens r'\ut+r\uo\tens [r\ut,r'\uo]\tens
r'\ut+r\uo\tens r'\uo\tens [r\ut,r'\ut]=0}
and having $\ad$-invariant symmetric part $2r_+=r+\tau(r)$, where $\tau$ is
transposition. We use the conventions and notation similar to \cite[Ch.
8]{Ma:book}, with $r=r\uo\tens r\ut$ denoting an element of $g\tens g$
(summation understood) and $r'$ denoting another distinct copy of $r$. We also
use $\delta\xi=\xi\o\tens\xi\t$ to denote the output in $\cg\tens\cg$ for
$\xi\in\cg$ (summation understood). A quasitriangular Lie bialgebra is called
factorisable if $2r_+$ is surjective when viewed as a map $\cg^*\to\cg$.

In view of the discussion above, we are interested in Lie-algebraic objects
living in the category ${}_{\cg}\CM$ of modules over a quasitriangular Lie
bialgebra $\cg$.
If $V$ is a $\cg$-module, we define its {\em infinitesimal braiding} to be the
operator
\[ \psi:V\tens V\to V\tens V, \quad \psi(v\tens w)=2r_+\la(v\tens w-w\tens v)\]
where $\la$ denotes the left action of $\cg$.

\begin{lemma} Let $\cb\in{}_{\cg}\CM$ be a $\cg$-covariant Lie algebra. Then
the associated $\psi:\cb\tens\cb\to \cb\tens\cb$ is a 2-cocycle in $\psi\in
Z^2_{\ad}(\cb,\cb\tens\cb)$.
\end{lemma}
\proof The proof that $\extd\psi=0$ is a straightforward computation in Lie
algebra cohomology. We use covariance of $\cb$ in the form: $xi\la[x,y]=[\xi\la
x,y]+[x,\xi\la y]$ for all $\xi\in \cg$. Then,
\align{(\extd\psi)(x,y,z)\equad&&=-\psi([x,y],z)+\psi([x,z],y)
-\psi([y,z],x)+\ad_x\psi(y,z)-\ad_y\psi(x,z)+\ad_z\psi(x,y)\\
&&=2r_+\la(-[x,y]\tens z+z\tens [x,y])+[x,2r_+\uo\la y]\tens r_+\ut\la
z+2r_+\uo\la y\tens [x,r_+\ut\la z]\\
&&\quad -[x,2r_+\uo\la z]\tens r_+\ut\la y+2r_+\uo\la z\tens [x,r_+\ut\la
y]+{\rm cyclic}\\
&&=-[2r_+\uo\la x,y]\tens z +2r_+\uo\la y\tens [x,r_+\ut\la z]-[x,2r_+\uo\la
z]\tens r_+\ut\la y\\
&&\quad + 2r_+\uo\la z\tens [r_+\ut\la x,y]+ {\rm cyclic}=0}
on using the cyclic invariance in $x,y,z$ and antisymmetry of the Lie bracket.
Note that this works for any element $2r_+\in \cg\tens\cg$ in the definition of
$\psi$. \endproof

\begin{defin} A {\em braided-Lie bialgebra} $\cb\in{}_{\cg}\CM$ is a
$\cg$-covariant
Lie algebra and $\cg$-covariant Lie coalgebra with cobracket
$\und\delta:\cb\to \cb\tens \cb$ obeying $\forall x,y\in\cb$,
\[\und\delta([x,y])=\ad_x\delta y-\ad_y\delta x-\psi(x\tens y);\quad
\psi=2r_+\la(\id-\tau),\]
i.e., $\und\delta$ obeys the coJacobi identity and $\extd\und\delta=\psi$.
\end{defin}

The definition is motivated from that of a braided group, where the coproduct
fails to be multiplicative up to a braiding $\Psi$\cite{Ma:bg}. The results in
the next section serve to justify it further.

\section{Lie versions of braided group theorems}

The existence of nontrivial quasitriangular Lie bialgebra structures is
known\cite{Dri} for all simple $\cg$ at least over $\C$. Our first theorem
ensures likewise the existence of braided-Lie bialgebras.

\begin{theorem} Let $i:\cg\to\cf$ be a map of Lie bialgebras, with $\cg$
quasitriangular. There is a braided-Lie bialgebra $\cb(\cg,\cf)$, the {\em
transmutation}
of $\cf$, living in ${}_{\cg}\CM$  with Lie algebra
$\cf$ and for all $x\in \cf$, $\xi\in\cg$,
\[ \und\delta x=\delta x + r\uo\la x\tens i(r\ut)-i(r\ut)\tens r\uo\la x,\quad
\xi\la x=[i(\xi),x].\]
\end{theorem}
\proof We first verify that $\und\delta$ as stated is indeed a $\cg$-module
map. Thus
\align{\und\delta(\xi\la x)\equad&&=\delta[i(\xi),x]+r\uo\xi\la x\tens
i(r\ut)-i(r\ut)\tens r\uo\xi\la x\\
&&=\xi\la\delta x-[x,i(\xi\o)]\tens
i(\xi\t)-i(\xi\o)\tens[x,i(\xi\t)]+r\uo\xi\la x\tens i(r\ut)-i(r\ut)\tens
r\uo\xi\la x\\
&&=\xi\la\delta x-[x,[i(\xi),i(r\uo)]]\tens i(r\ut)-[x,i(r\uo)]\tens
[i(\xi),i(r\ut)]\\
&&\quad -[i(\xi),i(r\uo)]\tens [x,i(r\ut)]-i(r\uo)\tens
[x,[i(\xi),i(r\ut)]]+r\uo\xi\la x\tens i(r\ut)-i(r\ut)\tens r\uo\xi\la x\\
&&=\xi\la\delta x+[\xi,r\uo]\la x\tens i(r\ut)-[x,i(r\uo)]\tens
[i(\xi),i(r\ut)]\\
&&\quad -[i(\xi),i(r\uo)]\tens [x,i(r\ut)]+i(r\uo)\tens [\xi,r\ut]\la
x+r\uo\xi\la x\tens i(r\ut)-i(r\ut)\tens r\uo\xi\la x\\
&&=\xi\la\delta x+r\uo\la x\tens [i(\xi),i(r\ut)]+[i(\xi),i(r\uo)]\tens r\ut\la
x+\xi r\uo\la x\tens i(r\ut)-i(r\ut)\tens \xi r\uo\la x\\
&&=\xi\la\und\delta x}
where we used the definitions of $\la$ and $\und\delta$ and the fact that $\cg$
is quasitriangular, so that $\delta\xi=[\xi,r\uo]\tens
r\ut+r\uo\tens[\xi,r\ut]$.

Antisymmetry of the output of $\und\delta$ is clear. Next we verify the
coJacobi identity,
\align{(\id\tens\und\delta)\und\delta x\equad&&+{\rm
cyclic}=(\id\tens\und\delta)\delta+ r\uo\la x\tens\und\delta
i(r\ut)-i(r\ut)\tens\und\delta(r\uo\la x)+{\rm cyclic}\\
&&\equad=x\o\tens r\uo\la x\t\tens i(r\ut)-x\o\tens i(r\ut)\tens r\uo\la
x\t+r\uo\la x\tens \delta i(r\ut)\\
&&+r\uo\la x\tens i([r'\uo,r\ut])\tens i(r'\ut)-r\uo\la x\tens i(r'\ut)\tens
i([r'\uo,r\ut])-i(r\ut)\tens r\uo\la x\o\tens x\t\\
&&-i(r\ut)\tens x\o\tens r\uo\la x\t-i(r\ut)\tens r\uo r'\uo\la x\tens
i(r'\ut)-i(r\ut)\tens r'\uo\la x\tens i([r\uo,r'\ut])\\
&&+i(r\ut)\tens i([r\uo,r'\ut])\tens r'\uo\la x+ i(r\ut)\tens i(r'\ut)\tens
r\uo r'\uo\la x+{\rm cyclic}}
using the definition of $\und\delta$ and the previous covariance result.
Several of the resulting terms cancel immediately. Using the quasitriangular
form of $\delta$ on $r\ut$ and the further freedom to cyclically rotate all
tensor products so that $x$ appears in the first factor, our expression becomes
\align{&&\equad =r\uo\la x\tens i([r\ut,r\uo])\tens i(r'\ut)+r\uo\la x\tens
i(r'\uo)\tens i([r\ut,r'\ut])\\
&&+r\uo\la x\tens i([r'\uo,r\ut])\tens i(r'\ut)-r\uo\la x\tens i(r\ut)\tens
i([r'\uo,r\ut])+[r\uo,r'\uo]\la x\tens i(r\ut)\tens i(r'\ut)\\
&&+r'\uo\la x\tens i(r\ut)\tens i([r\uo,r'\ut])-r'\uo\la x\tens
i([r\uo,r'\ut])\tens i(r\ut)+{\rm cyclic}\\
&&\equad=( (\ )\la x\tens i\tens i)\left([r\uo,r'\uo]\tens r\ut\tens
r'\ut+r\uo\tens [r\ut,r\uo]\tens r\ut+r\uo\tens r'\uo\tens
[r\ut,r'\ut]\right)+{\rm cyclic}\\
&&\equad=0}
by the CYBE (\ref{CYBE}).

Finally, we prove that $\extd\und\delta=\psi$. Thus,
\align{&&\equad\und\delta([x,y])=\delta[x,y]+r\uo\la[x,y]\tens
i(r\ut)-i(r\ut)\tens r\uo\la[x,y]\\
&&=\ad_x\delta y-\ad_y\delta x + [r\uo\la x,y]\tens i(r\ut)+[x,r\uo\la y]\tens
i(r\ut)-i(r\ut)\tens[r\uo\la x,y]-i(r\ut)\tens [x,r\uo\la y]\\
&&=\ad_x\und\delta y-\ad_y\und\delta x-r\uo\la y\tens [x,i(r\ut)]+r\uo\la
x\tens [y,i(r\ut)]+[x,i(r\ut)]\tens r\uo\la y - [y,i(r\ut)]\tens r\uo\la x\\
&&=\ad_x\und\delta y -\ad_y\und\delta x-2r_+\la(x\tens y-y\tens x)}
as required. We used the definitions of $\und\delta$ and $\la$. \endproof

\begin{corol}
Every quasitriangular Lie bialgebra $\cg$ has a braided version
$\und\cg\in{}_{\cg}\CM$ by $\ad$,
the same Lie bracket, and
\eqn{kk}{\und\delta x= 2r_+\uo\tens [x,r_+\ut].}
\end{corol}
\proof We take the identity map $i=\id:\cg\to\cg$ and $\und\cg=\cb(g,g)$. Its
braided-Lie cobracket from Theorem~3.1 is $\und\delta x=[x,r\uo]\tens
r\ut+r\uo\tens[x,r\ut]+r\uo\la x\tens r\ut-r\ut\tens r\uo\la x$ using the
quasitriangular form of $\delta$. \endproof

The corollary ensures the existence of non-trivial braided-Lie bialgebras since
nontrivial quasitriangular Lie bialgebras are certainly known.

\begin{example} Let $\cg$ be a finite-dimensional factorisable Lie bialgebra.
Then $\und\delta$ in Corollary~3.2 is equivalent under the isomorphism
$2r_+:\cg^*\isom\cg$ to the Kirillov-Kostant Lie cobracket on $\cg^*$ (defined
as the dualisation of the Lie bracket $\cg\tens\cg\to\cg$). The braided-Lie
bialgebra $\und\cg$ is self-dual.
\end{example}
\proof It is well known that for any Lie algebra the vector space $\cg^*$
acquires a natural Poisson bracket structure. Considering $\cg$ as a subset of
the functions on $g^*$, this Kirillov-Kostant Poisson bracket is
$\{\xi,\eta\}(\phi)=\<\phi,[\xi,\eta]\>$ where $\<\ ,\ \>$ denotes evaluation
and $\xi,\eta\in\cg$, $\phi\in\cg^*$. The associated Lie coalgebra structure
$\und\delta :\cg^*\to\cg^*\tens\cg^*$ is defined by
$\{\xi,\eta\}(\phi)=\<\xi\tens\eta,\und\delta \phi\>$ and is therefore the
dualisation of the Lie bracket of $\cg$. We call it the Kirillov-Kostant Lie
coalgebra structure on $\cg^*$.

Let $K(\phi)=2r_+\uo\<\phi,r_+\ut\>$ denote the isomorphism $K:\cg^*\isom\cg$
resulting from our factorisability assumption.
Then
\align{\<\xi\tens\eta,(K^{-1}\tens K^{-1})\und\delta K(\phi)\>\equad
&&=\<\xi,K^{-1}(2r_+\uo)\>\<\eta,K^{-1}([K(\phi),r_+\ut])\>\\
&&=\<K^{-1}(\xi),2r_+\uo\>\<[\eta,K(\phi)],K^{-1}(r_+\ut)\>\\
&&=\<[\eta,K(\phi)],K^{-1}(\xi)\>=\<K(\phi),K^{-1}([\xi,\eta])\>
=\<\phi,[\xi,\eta]\>.}
We used symmetry and ad-invariance of $K$ as an element of $g\tens g$, with its
corresponding property
$\<\eta,K^{-1}([\xi,\zeta])\>=\<[\eta,\xi],K^{-1}(\zeta)\>$  $\forall
\xi,\eta,\zeta\in\cg$, for the map $K:\cg^*\to\cg$.

Next, we give $\cg^*$ with the above Kirillov-Kostant Lie cobracket
$\delta\phi=\phi\o\tens\phi\t$ (dual to the Lie algebra of $\cg$) a Lie bracket
and $\cg$-module structure
\eqn{g*alg}{[\phi,\chi]=\chi\o 2r_+(\phi,\chi\t),\quad
\xi\la\phi=\phi\o\<\phi\t,\xi\>}
for all $\xi\in\cg$, $\phi,\chi\in\cg^*$ and with $2r^+$ viewed as a map
$\cg^*\tens\cg^*\to k$. Then $\cg^*$ becomes a braided-Lie bialgebra in
${}_{\cg}\CM$, which we denote $\und\cg^*$. Its Lie cobracket is
$\und\delta=\delta$ the dual of the Lie bracket of $\und\cg$ (since this is the
same as that of $\cg$), and its Lie bracket is dual to the Lie cobracket of
$\cg$ in Corollary~3.2 since
\[\<\xi,[\phi,\chi]\>=\<\xi,\chi\o\>\<2r_+(\phi),\chi\t\>
=\<[\xi,2r_+(\phi)],\chi\>=\<2r_+\ut\tens[\xi,r_+\uo],\phi\tens\chi\>
=\<\und\delta\xi,\phi\tens\chi\>\]
for all $\xi\in\cg$ and $\phi,\chi\in\cg^*$. On the other hand,
\[
\<\xi,[\phi,\chi]\>=\<\xi,\chi\o\>\<\chi\t,K(\phi)\>=\<[\xi,K(\phi),\chi\>
=\<\xi,K^{-1}([K(\phi),K(\chi)])\>\]
hence $2r_+:\und\cg^*\to\und\cg$ is an isomorphism of braided-Lie bialgebras.
\endproof

This is the Lie analogue of the theorem that braided groups obtained by full
transmutation of factorisable quantum groups are self-dual via the quantum
Killing form\cite{Ma:skl}. Also, the fact that the data corresponding to the
original Lie cobracket on $\cg$ does not enter into $\und\cg$ corresponds in
braided group theory to transmutation rendering a quasitriangular Hopf algebra
braided-cocommutative. There is also a theory of quasitriangular braided-Lie
bialgebras of which the more general $\cb(\cg,\cf)$ are examples when $\cf$ is
itself quasitriangular. The braided-quasitriangular structure is the difference
of the quasitriangular structures on $\cf,\cg$ as the Lie version of results in
\cite{Ma:tra}.

For use later on, the general duality for braided-Lie bialgebras relevant to
Example~3.3 is given by
\begin{lemma} If $\cb\in{}_{\cg}\CM$ is a finite-dimensional braided-Lie
bialgebra then there is  a dual braided-Lie bialgebra $\cb^*\in{}_{\cg}\CM$. It
is built on the vector space $\cb^*$ with action $(\xi\la\phi)(x)=-\phi(\xi\la
x)$ and Lie (co)bracket structure maps given by dualisation.
\end{lemma}
\proof The Jacobi and coJacobi (and antisymmetry) axioms are clear by
dualisation, as is the specified left action on $\cb^*$.  The induced
infinitesimal braiding on the dual is the usual dual:
\[ \<x\tens y,\psi_{\cb^*}(\phi,\chi)\>=\<x\tens
y,2r_+\la(\phi\tens\chi-\chi\tens\phi)\>=\<\psi(x\tens y),\phi\tens\chi\>\]
for all $x,y\in\cb$ and $\phi,\chi\in\cb^*$. Moreover, the map
$\extd\und\delta$ for $\cb$ dualises to $\extd\und\delta$ for $\cb^*$. The
proof is identical to the proof  that the dual of a usual Lie bialgebra is a
Lie bialgebra (see \cite{Ma:book} for details). Hence $\extd\und\delta=\psi$
for $\cb^*$ by dualisation of this relation for $\cb$. \endproof

Note also that if $\cb\in{}_{\cg}\CM$ is any braided-Lie bialgebra then so is
$\cb^{\rm op/cop}$ with opposite bracket and cobracket, in the same category.
This is because the covariance conditions on the Lie bracket and cobracket are
each linear in those structures and hence valid even with the additional minus
signs in either case. Meanwhile, $\extd\und\delta$ is linear in $\und\delta$
and linear in the Lie bracket, hence invariant when both are changed by a minus
sign. The infinitesimal braiding does not involved either the bracket or
cobracket and is invariant. Applying this observation to $\cb^*$ gives us
another dual, $\cb^\star$. This is the Lie analogue of the more
braided-categorical which is more natural in the theory of braided groups. In
the Lie setting, however, we have $\cb^\star\isom\cb^*$ by $x\mapsto -x$, so
can work entirely with $\cb^*$. We also conclude, in passing, that $\cb^{\rm
op}$ and $\cb^{\rm cop}$ are braided-Lie bialgebras in the category of modules
over the opposite quasitriangular Lie bialgebra (i.e. with quasitriangular
structure $-r_{21}$ in place of $r$).

We consider now the adjoint direction to Theorem~3.1, to associate to a
braided-Lie bialgebra an ordinary Lie bialgebra. The quantum group version of
this result\cite{Ma:bos} has been used to construct inhomogeneous quantum
groups\cite{Ma:poi}.

\begin{theorem}
Let $\cb\in{}_{\cg}\CM$ be a braided-Lie bialgebra. Its {\em
bosonisation} is the Lie bialgebra $\cb\lbiprod \cg$ with $\cg$ as sub-Lie
bialgebra, $\cb$ as sub-Lie algebra and
\eqn{liebos}{[\xi,x]=\xi\la x,\quad \delta x=\und\delta x+r\ut\tens r\uo\la
x-r\uo\la x\tens r\ut,\quad\forall \xi\in\cg,\ x\in\cb.}
\end{theorem}
\proof The Lie algebra structure of $\cb\lbiprod\cg$ is constructed as a
semidirect sum by the given action of $\cg$ on $\cb$. The coassociativity of
the Lie cobracket may be verified directly from the CYBE along the lines of the
proof of coassociativity in Theorem~3.1.  The line of deduction is reversed by
the formulae are similar. That the result is a Lie bialgebra has three cases.
For $\xi,\eta\in\cg$ we have $\delta([\xi,\eta])$ as required since  $\cg$ is a
Lie bialgebra. For the mixed case we have
\align{\delta([\xi,x])\equad&&=\delta(\xi\la x)=\und\delta(\xi\la x)+r\ut\tens
r\uo\xi\la x-r\uo\xi\la x\tens r\ut\\ -\ad_\xi\delta
x\equad&&=-\xi\la\und\delta x-[\xi,r\ut]\tens r\uo\la x -r\ut\tens\xi r\uo\la
x+\xi r\uo\la x\tens r\ut+r\uo\la x\tens[\xi,r\ut]\\
\ad_x\delta\xi\equad&&=-\xi\o\la x\tens\xi\t-\xi\o\tens\xi\t\la x\\
&&=[r\uo,\xi]\la x\tens r\ut+r\uo\la x\tens[r\ut,\xi]-[r\uo,\xi]\tens r\ut\la
x-r\uo\tens[r\ut,\xi]\la x.}
We used the definitions and, in the last line, the form of $\delta\xi$ as a
quasitriangular Lie bialgebra. Adding these expressions, we obtain
\[ \delta([\xi,x])-\ad_\xi\delta
x+\ad_x\delta\xi=\und\delta\xi-\xi\la\und\delta\xi +[2r_+\ut,\xi]\tens
r_+\uo\la x+2r_+\ut\tens[r_+\uo,\xi]\la x=0\]
by covariance of $\und\delta$ and $\ad$-invariance of $2r_+$. The remaining
case is
\align{\delta([x,y])\equad&&=\und\delta([x,y])+r\ut\tens
r\uo\la[x,y]-r\uo\la[x,y]\tens r\ut\\
&&=(\ad_x\und\delta y+r\ut\tens \ad_x(r\uo\la y) -\ad_x(r\uo\la y)\tens
r\ut-(x\swap y))-\psi(x\tens y)\\
&&=(\ad_x\delta y -\ad_x(r\ut)\tens r\uo\la y+r\uo\la y\tens\ad_x(r\ut)
-(x\swap y))-\psi(x\tens y)\\
&&=\ad_x\delta y-\ad_y\delta x}
on writing $\ad_x(r\ut)=-r\ut\la x$ and comparing with the definition of
$\psi$. We used the braided-Lie bialgebra property of $\und\delta$. \endproof

The construction in the bosonisation theorem can also be viewed as a special
case of a more general construction for Lie bialgebras which are semidirect
sums as Lie algebras and Lie coalgebras by a simultaneous Lie action and Lie
coaction. We call such Lie algebras {\em bisum} Lie algebras. They are the
analogue of {\em biproduct} Hopf algebras in \cite{Rad:str}. In the general
case one only needs covariance under a Lie bialgebra, not necessarily
quasitriangular. However, any Lie bialgebra has a Drinfeld double\cite{Dri:ham}
which is quasitriangular. In order to explain these topics we need quite a bit
more formalism. Firstly, if $\cf$ is any Lie coalgebra, we have a notion of
left {\em Lie coaction} on a vector space $V$. This is a map $\beta:V\to
\cf\tens V$ such that\cite{Ma:book}
\eqn{coact}{ (\delta\tens\id)\circ\beta=((\id-\tau)\tens\id)\circ(\id\tens
\beta)\circ\beta.}
The category of left Lie comodules is denoted ${}^{\cg}\CM$ and is monoidal in
the obvious derivation-like way. Morphisms are defined as linear maps
intertwining the Lie coactions, again in the obvious way.

\begin{lemma} Let $\cf$ be a Lie bialgebra. There is a monoidal category of
{\em crossed modules} ${}^{\cf}_{\cf}\CM$  having as objects vector spaces $V$
which are
are simultaneously $\cf$-modules $\la:\cf\tens V\to\tens V$ and $\cf$-comodules
$\beta:V\to\cf\tens V$ obeying $\forall \xi\in\cf, v\in V$,
\[ \beta(\xi\la v)=([\xi, ]\tens\id+\id\tens \xi\la)\beta(v)+(\delta \xi)\la
v.\]
It can be identified when $\cf$ is finite-dimensional with the category
${}_{D(\cf)}\CM$ where $D(\cf)$ is the Drinfeld double\cite{Dri:ham}. Writing
$\beta(v)=v\bo\tens v\bt$, the corresponding infinitesimal braiding on any
object  $V\in {}_{\cf}^{\cf}\CM$ is
\[\psi(v\tens w)= w\bo\la v\tens w\bt-v\bo\la w\tens v\bt- w\bt\tens w\bo\la v+
v\bt\tens v\bo\la w.\]
\end{lemma}
\proof Morphisms in ${}^{\cf}_{\cf}\CM$ are maps intertwining both the Lie
action and the Lie coaction. We start with $\cf$ finite-dimensional and use
Drinfeld's formulae for $D(\cf)$ in the conventions in \cite{Ma:book}, which
contains the Lie algebras $\cf$ and $\cf^{*\rm op}$ with the cross relations
$[\xi,\phi]=\phi\o\<\phi\t,\xi\>+\xi\o\<\phi,\xi\t\>$ for all  $\xi\in\cf$ and
$\phi\in\cf^*$.
A left module of $D(\cf)$ therefore means a vector space which is a left
$\cf$-module and a right $\cf^*$-module, obeying $\xi\la(v\ra\phi)-(\xi\la
v)\ra\phi=v\ra\phi\o\<\phi\t,\xi\>+\xi\o\la v\<\phi,\xi\t\>$.  Next we view the
right action of $\cf^*$ as, equivalently, a left coaction of $\cf$ by
$v\ra\phi=\<\phi,v\bo\>v\bt$. Inserting this, we have the condition
\[ \xi\la v\bt \<\phi,v\bo\>-\<\phi,(\xi\la v)\bo\>(\xi\la
v)\bt=\<\phi,[v\bo,\xi]\>v\bt+\xi\o\la v \<\phi,\xi\t\>\]
for all $\phi$. We wrote the Lie cobracket of $\cf^*$ in terms of the Lie
algebra of $\cf$ here. This is the condition stated for $\beta$, which
manifestly makes sense even for infinite-dimensional Lie algebras. It is easy
to check that the category is well defined an monoidal even in this case. In
the same spirit, $D(\cf)$ has a quasitriangular structure given by the
canonical element for the duality pairing\cite{Dri:ham}. Then
$2r_+=\sum_af^a\tens e_a+e_a\tens f^a$
where $\{e_a\}$ is a basis of $\cf$ and $\{f^a\}$ is a dual basis. Hence the
infinitesimal braiding in ${}_{D(\cf)}\CM$ is
\[ \psi(v\tens w)=2r_+\la(v\tens w-w\tens v)=\<f^a,v\bo\>v\bt\tens e_a\la
w+e_a\la v\tens\<f^a,w\bo\>w\bt -(v\swap w)\]
when the left action of $\cf^{*\rm op}$ is of the form given by a left coaction
of $\cf$.
This gives $\psi$ as stated. Note that the use of a coaction to reformulate an
action of the dual in the infinite dimensional case is a completely routine
procedure in Hopf algebra theory; we have given the details here since the Lie
version is less standard; the category of Lie crossed modules
${}^{\cf}_{\cf}\CM$ should not be viewed as anything other than a {\em version}
of the ideas behind Drinfeld's double construction. \endproof

The resulting map $\psi$ is well-defined even in the infinite dimensional case;
we call it the  infinitesimal braiding of the category ${}_{\cf}^{\cf}\CM$ of
crossed $\cf$-modules and define a braided-Lie bialgebra in ${}_{\cf}^{\cf}\CM$
with respect to this.

\begin{theorem} Let $\cf$ be a Lie bialgebra and let $\cb\in{}_{\cf}^{\cf}\CM$
be a braided-Lie bialgebra. The {\em bisum} Lie
bialgebra $\cb\lbiprod \cf$ has semidirect Lie bracket/cobracket and projects
onto $\cf$. Conversely, any Lie bialgebra $\cg$ with a split Lie bialgebra
projection $\cg{ {\pi\atop\to}\atop {\hookleftarrow\atop i}}\cf$ is of this
form, with $\cb=\ker\pi$ and braided-Lie bialgebra structure given by
$\cb\subset\cg$ as a Lie algebra and
\[ \xi\la=\ad_{i(\xi)},\quad \beta=(\pi\tens\id)\circ\delta,\quad
\und\delta=(\id-i\circ\pi)^{\tens 2}\circ\delta.\]
\end{theorem}
\proof  In the forward direction, since $\cb$ is covariant under an action of
$\cf$ we can make, as usual, a semidirect sum $\cb\lcross\cf$. The bracket on
general elements of the direct sum vector space is
$[x\oplus\xi,y\oplus\eta]=([x,y]+\xi\la y-\eta\la x)\oplus [\xi,\eta]$ as
usual. On the other hand, since the Lie coalgebra of $\cb$ is covariant under a
Lie coaction of $\cf$, one may make a semidirect Lie coalgebra
$\cb\lcocross\cf$ with\cite{Ma:book}
\eqn{lcocross}{ \delta(x\oplus\xi)=\delta\xi+\und\delta
x+(\id-\tau)\circ\beta(x)}
where $\und\delta$ is the Lie cobracket of $\cb$. The required covariance of
the Lie coalgebra under the coaction here is
\eqn{comodcoalg}
{(\id\tens\und\delta)\circ\beta=(\id\tens(\id-\tau))
\circ(\beta\tens\id)\circ\und\delta}
and ensures that $\delta$ on $\cb\lcocross\cf$ obeys the coJacobi identity.

The further covariance assumptions on $\cb$ are that its Lie bracket is
covariant under the Lie coaction and its Lie cobracket is covariant under the
Lie action.  These assumptions are all needed to show that the semidirect Lie
bracket and Lie cobracket form a Lie bialgebra $\cb\lbiprod\cf$. The case
$\delta([\xi,\eta])$ is clear since $\cf$ is a Lie subalgebra and Lie
subcoalgebra. The mixed case is
\align{\delta([\xi,x])\equad&&=\delta(\xi\la x)=\und\delta(\xi\la x)+(\xi\la
x)\bo\tens (\xi\la x)\bt-(\xi\la x)\bt\tens (\xi\la x)\bo\\
&&=\xi\la\und\delta x+\ad_\xi((\id-\tau)\circ\beta(x))+\xi\o\tens\xi\t\la
x-\xi\t\tens\xi\o\la x\\
&&=\ad_\xi\delta x-\ad_x\delta\xi}
as required. We used the definition of $\delta$ for $\cb\lcocross\cf$, the
covariance of $\und\delta$ under $\xi$, antisymmetry of $\delta\xi$ and
$\ad_\xi(x)=\xi\la x=-\ad_x(\xi)$ when viewed in the Lie algebra
$\cb\lcross\cf$. The remaining case is
\align{\delta([x,y])\equad &&=\und\delta([x,y])+(\id-\tau)\beta([x,y])\\
&&=\ad_x\und\delta y+y\bt\tens y\bo\la x-y\bo\la x\tens y\bt+x\bo\tens
[x\bt,y]-[x\bt,y]\tens x\bo -(x\swap y)\\
&&=\ad_x\und\delta y-y\bo\la x\tens y\bt+y\bo\tens[x,y\bt]-[x,y\bt]\tens
y\bo+y\bt\tens y\bo\la x-(x\swap y)\\
&&=\ad_x\und\delta y+[x,y\bo]\tens y\bt+y\bo\tens[x,y\bt]-[x,y\bt]\tens
y\bo-y\bt\tens[x,y\bo]-(x\swap y)\\
&&=\ad_x\delta y-\ad_y\delta x}
where we used  $\beta([x,y])=y\bo\tens [x,y\bt]-x\bo\tens[y,x\bt]$ (covariance
of the Lie bracket of $\cb$ under the Lie coaction) and the assumption that
$\und\delta$ is a braided-Lie bialgebra in ${}^{\cf}_{\cf}\CM$ with
infinitesimal braiding from Lemma~3.6.
We then used the crossed module compatibility condition also from Lemma~3.6 and
$\xi\la x=-\ad_x(\xi)$ to recognise the required result. It is easy to see that
the projection $\cb\lbiprod\cf\to\cf$ defined by setting elements of $\cb$ to
zero is a Lie bialgebra map covering the inclusion $\cf\subset \cb\lbiprod\cf$.

In the converse direction, we assume a split projection, i.e. a surjection
$\pi:\cg\to\cf$ between Lie bialgebras covering an inclusion $i:\cf\to \cg$ of
Lie bialgebras (so that $\pi\circ i=\id$). We define $\cb=\ker\pi$. Since this
is a Lie ideal, it both forms a sub-Lie algebra of $\cg$ and is covariant under
the action of $\cf$ given by pull-back along $i$ of $\ad$. Moreover, $\cg$
coacts on itself by its Lie cobracket $\delta$ (the adjoint coaction of any Lie
bialgebra on itself) and hence push-out along $\pi$ is an $\cf$-coaction
$\beta=(\pi\tens\id)\circ\delta$, which restricts to $\cb$ since
$(\id\tens\pi)\beta(x)=(\pi\tens\pi)\delta x=\delta\pi(x)=0$ for $x\in\ker\pi$.
This Lie action and Lie coaction fit together to form a Lie crossed module,
\align{\beta(\xi\la x)\equad
&&=(\pi\tens\id)\circ\delta([i(\xi),x])=(\pi\tens\id)(\ad_{i(\xi)}\delta
x-\ad_x\delta i(\xi))\\
&&=\pi([i(\xi),x\o])\tens x\t+\pi(x\o)\tens[i(\xi),x\t]-\pi([x,i(\xi)\o])\tens
i(\xi)\t-\pi(i(\xi)\o)\tens[x,i(\xi)\t]\\
&&=[x,\pi(x\o)]\tens x\t+\pi(x\o)\tens\xi\la x\t+\xi\o\tens\xi\t\la x\\
&&=([\xi,\ ]\tens\id+\id\tens\xi\la)\beta(x)+(\delta\xi)\la x}
as required. We used that $i$ is a Lie coalgebra map and $\pi$ a Lie algebra
map, along with
$x\in\ker\pi$ to kill the term with $\pi([x,i(\xi)\o])$.

Finally, we give $\cb$ a Lie cobracket $\und\delta$ as stated. Writing
$p=\id-i\circ\pi$, we have \align{(\id\tens\und\delta)\circ\und\delta x\equad
&&=p(x\o)\tens p(p(x\t)\o)\tens p(p(x\t)\t)\\
&&=(p\tens p\tens p)(\id\tens\delta)\delta x-p(x\o)\tens
p(i\circ\pi(x\t)\o)\tens p(i\circ\pi(x\t)\t)\\
&&=(p\tens p\tens p)(\id\tens\delta)\delta x}
since  $i\circ\pi$ is a Lie coalgebra map and $p\circ i\circ\pi=0$. Hence
$\und\delta$ obeys the coJacobi identity since $\delta$ does. Moreover, for all
$x,y\in\ker\pi$,
\align{\und\delta([x,y])\equad&&=(\id-i\circ\pi)\tens(\id-i\circ\pi)\delta x\\
&&=[x,y\o]\tens y\t+y\o\tens[x,y\t]-[x,y\o]\tens
i\circ\pi(y\t)-i\circ\pi(y\o)\tens[x,y\t]}
since $i\circ\pi([x,y\t])=0$ etc., as $i\circ\pi$ is a Lie algebra maps. Also,
from Lemma~3.6 and the form of $\beta$ and antisymmetry of $\delta$ we have
\align{\psi(x\tens y)\equad&&=[i\circ\pi(y\t),x]\tens
y\t+y\o\tens[i\circ\pi(y\t),x]-(x\swap y).}
Then,
\align{\ad_x\und\delta y\equad &&-\ad_y\und\delta
x=[x,(\id-i\circ\pi)(y\o)]\tens (\id-i\circ\pi)(y\t)-(x\swap y)\\
&&=[x,y\o]\tens y\t+y\o\tens [x,y\t]-[x,i\circ\pi(y\o)]\tens y\t-[x,y\o]\tens
i\circ\pi(y\t)\\
&&\quad -i\circ\pi(y\o)\tens [x,y\t]-y\o\tens[x,i\circ\pi(y\t)]-(x\swap
y)=\psi(x\tens y)+\und\delta([x,y])}
as required. The additional terms $i\circ\pi(y\o)\tens [x,i\circ\pi(y\t)]$ etc.
vanish as $i\circ\pi$ is a Lie coalgebra map and $x,y\in\ker\pi$. Hence
$\cb=\ker\pi$ becomes a braided-Lie bialgebra in ${}_{\cf}^{\cf}\CM$. One may
then verify that the bisum Lie bialgebra $\cb\lbiprod\cf$ coincides with  $\cg$
viewed as a direct sum $\cb\oplus\cf$ of vector spaces according to the
projection $i\circ\pi$. \endproof

This is the Lie analogue of the braided groups interpretation\cite{Ma:skl} of
Radford's theorem\cite{Rad:str}. It tells us that braided-Lie bialgebras are
rather common as they arise whenever we have a projection of ordinary Lie
bialgebras. Finally, we provide the Lie analogue of the functor\cite{Ma:dou}
which connects biproducts and bosonisation.

\begin{lemma} Let $\cg$ be a quasitriangular Lie bialgebra. There is a monoidal
functor
${}_{\cg}\CM\to{}^{\cg}_{\cg}\CM$ respecting the infinitesimal braidings. It
sends an action $\la$ to a pair $(\la,\beta)$ where $\beta=r_{21}\la$, the {\em
induced Lie coaction}.  The bosonisation of $\cb\in{}_{\cg}\CM$ in Theorem~3.5
can thereby be viewed as an example of a biproduct in Theorem~3.7.
\end{lemma}
\proof We first verify that $\beta(v)=r\ut\tens r\uo\la v$ defines a Lie
coaction for any $\cg$-module $V\owns v$. This follows immediately from the
identity $(\id\tens\delta)r=[r\uo,r'\uo]\tens r'\ut\tens r\ut$ holding for any
quasitriangular Lie bialgebra (following from the CYBE and $\delta=\extd r$).
Thus,
\align{(\id\tens\delta)\beta(v)\equad &&= r\ut\tens r'\ut\tens r'\uo r\uo\la v-
r'\ut\tens r\ut\tens r'\uo r\uo\la v\\
&&=((\id-\tau)\tens\id)\circ(\id\tens \beta)\circ\beta(v)}
as required. This fits together with the given action to form a Lie crossed
module as
\align{\beta(\xi\la v)\equad&&=r\ut\tens r\uo\xi\la v=r\ut\tens [r\uo,\xi]\la
v+r\ut\tens\xi r\uo\la v\\
&&=(\delta\xi)\la v+[\xi,r\ut]\tens r\uo\la v+(\id\tens\xi\la)\beta(v)}
as required, using the quasitriangular form of $\delta\xi$. More trivially, a
morphism $\phi:V\to W$ in ${}_{\cg}\CM$ is automatically an intertwiner of the
induced coactions (since $r\ut\tens r\uo\la\phi(v)=r\ut\tens \phi(r\uo\la v)$)
and hence a morphism in ${}_{\cg}^{\cg}\CM$. It us also clear that the functor
respects tensor products. In this way, ${}_{\cg}\CM$ is a full monoidal
subcategory.

Finally, we check that the infinitesimal braidings coincide. Computing $\psi$
from Lemma~3.6 in the image of the functor, we have
\[\psi(v\tens w)=r\ut\la v\tens r\uo\la w+r\uo\la v\tens r\ut\la w - (v\swap
w)=2r_+\la(v\tens w-w\tens v)\]
as required. From the form of the Lie cobracket in the bosonisation
construction, it is clear that it can be viewed as a semidirect Lie coalgebra
by the induced action, i.e it can be viewed as a nontrivial construction for
examples of bisum Lie algebras. \endproof

There is a dual theory of {\em dual quasitriangular} (or coquasitriangular) Lie
bialgebras\cite{Ma:book} where the Lie bracket has a special form
\eqn{dqua}{[\xi,\eta]=\xi\o r(\xi\t,\eta)+\eta\o
r(\xi,\eta\t),\quad\forall\xi,\eta\in\cg,}
defined by a dual quasitriangular structure $r:\cg\tens\cg\to k$. This is
required to obey the CYBE in a dual form
\eqn{dcybe}{r(\xi,\eta\o)r(\eta\t,\zeta)+r(\xi\o,\eta)r(\xi\t,\zeta)
+r(\xi,\zeta\o)r(\eta,\zeta\t)=0,\quad\forall \xi,\eta,\zeta\in\cg}
and $2r_+$ is required to be invariant under the adjoint Lie coaction
($=\delta$, the Lie cobracket) according to
$r_+(\xi,\eta\o)\eta\t+r_+(\xi\o\tens\eta)\xi\t=0$.  All of the above theory
goes through in this form. Thus, ${}^{\cg}\CM$ has, by definition, an
infinitesimal braiding defined by
\eqn{dpsi}{\psi(v\tens w)=r(v\bo,w\bo)(v\bt\tens w\bt-w\bt\tens v\bt)}
with respect to which we define a braided-Lie bialgebra in ${}^{\cg}\CM$. The
Lie comodule transmutation  theory associates to a map $\cf\to\cg$ of Lie
bialgebras with $\cg$ dual quasitriangular,   a braided-Lie bialgebra
$\cb(\cf,\cg)\in{}^{\cg}\CM$.

For example, the Lie comodule version of Corollary~3.2 is
$\und\cg\in{}^{\cg}\CM$ with the same Lie cobracket as $\cg$, the adjoint
coaction $\delta$ and
\eqn{dtrans}{ [\xi,\eta]=\eta\o 2r_+(\xi,\eta\t)\quad\forall\xi,\eta\in\cg.}
A concrete example is provided by $\cg^*$ when $\cg$ is finite-dimensional
quasitriangular. Then $\cg^*$ is dual quasitriangular and its transmutation
$\und{\cg^*}$ coincides with $(\und\cg)^*$ in (\ref{g*alg}) in Example~3.3.

Similarly, there is a functor ${}^{\cg}\CM\to {}^{\cg}_{\cg}\CM$ sending a Lie
coaction by $\cg$ to a crossed module with an induced action $\xi\la v=r(
v\bo,\xi)v\bt$ and respecting the infinitesimal braidings. A braided-Lie
bialgebra in $\cb\in{}^{\cg}\CM$ has a Lie bosonisation $\cb\lbiprod\cg$ given
by a semidirect Lie cobracket by the given Lie coaction and semidirect Lie
bracket given by the induced action. All of this dual theory follows rigorously
and automatically by writing all constructions in terms of equalities of linear
maps and then reversing all arrows. Such dualisation of theorems is completely
routine in the theory of Hopf algebras, and similarly here. Hence we do not
need to provide a separate proof of these assertions. Note that dualisation of
theorems should not be confused with the dualisation of given algebras and
coalgebras, which can be far from routine.

\begin{example} Let $\cg$ be a finite-dimensional quasitriangular Lie bialgebra
and $\und\cg^*$ the dual of its transmutation. Its bosonisation
$\und\cg^*\lbiprod\cg$ is isomorphic as a Lie bialgebra to the Drinfeld double
$D(\cg)$.
\end{example}
\proof The required isomorphism $\theta:D(\cg)\to \und\cg^*\lbiprod \cg$ is
$\theta(\phi)=\phi-r\ut\<\phi,r\uo\>$ and $\theta(\xi)=\xi$ for $\xi\in\cg$ and
$\phi\in\cg^*$. We check first that it is a Lie algebra map. The $[\xi,\eta]$
case is automatic as $\cg$ is a sub-Lie algebra on both sides. The mixed case
is
\align{[\theta(\xi),\theta(\phi)]_{\rm
bos}\equad&&=[\xi,\phi-r\ut\<\phi,r\uo\>]_{\rm
bos}=\xi\la\phi-[\xi,r\ut]\<\phi,r\uo\>\\
&&=\phi\o\<\phi\t,\xi\>-\xi\o\<\phi,\xi\t\>+r\ut\<\phi,[\xi,r\uo]\>\\
&&=\theta(\phi\o\<\phi\t,\xi\>+\xi\o\<\phi,\xi\t\>)=\theta([\xi,\phi])}
where $[\ ,\ ]_{\rm bos}$ is the Lie bracket of $\und\cg^*\lbiprod\cg$. We use
the definition of $\theta$, the quasitriangular form of $\delta\xi$, the action
$\xi\la\phi=\phi\o\<\phi\t,\xi\>$ for $\und\cg^*$ and the cross relations in
$D(\cg)$ (as recalled in Lemma~3.6)
to recognise the result. The remaining case is
\align{[\theta(\phi),\theta(\chi)]_{\rm
bos}\equad&&=[\phi-r\ut\<\phi,r\o\>,\chi-r'\ut\<\chi,r'\uo\>]\\
&&=[r\uo,r'\ut]\<\phi,r\uo\>\<\chi,r'\uo\>-r\ut\la\chi\<\phi,r\uo\>
+r\ut\la\phi\<\chi,r\uo\>+\chi\o\<2r_+,\phi\tens\chi\t\>\\
&&=[r\uo,r'\ut]\<\phi,r\uo\>\<\chi,r'\uo\>+\chi\o
\<r,\chi\t,\phi\>+\phi\o\<r,\chi\tens\phi\t\>\\
&&=[\chi,\phi]-r\ut\<\delta
r\uo,\chi\tens\phi\>=[\chi,\phi]-r\ut\<[\chi,\phi],r\uo\>=\theta([\chi,\phi])}
as required since $D(\cg)$ contains $\cg^{\rm op}$ as a sub-Lie algebra. We
used the definition of $\theta$ and the Lie bracket (\ref{g*alg}) of
$\und\cg^*$ as a sub-Lie algebra of the bosonisation.
We then used form of the action $r\uo\la\chi$ etc. and combined the result with
the $2r_+$ term to recognise the Lie bracket $[\chi,\phi]$ (as in (\ref{dqua}))
of the dual quasitriangular Lie bialgebra $\cg^*$. We also use the
quasitriangular form of $\cg$ to recognise $\delta r\uo$.

Next, we verify that $\theta$ is a Lie coalgebra map. This is automatic on
$\xi\in\cg$ as a sub-Lie bialgebra on both sides. The remaining case is
\align{\delta_{\rm bos}\theta(\phi)\equad &&=\delta_{\rm bos}\phi-\delta
r\ut\<\phi,r\uo\>\\
&&=\und\delta\phi+r\ut\tens r\uo\la\phi+r\uo\la\phi\tens r\ut-r\ut\tens
r'\ut\<\phi,[r\uo,r'\uo]\>\\
&&=\delta\phi+r\ut\tens\phi\o\<\phi\t,r\uo\>-\phi\o\<\phi\t,r\uo\>\tens
r\ut-r\ut\tens r'\ut\<\phi,[r\uo,r'\uo]\>\\
&&=(\phi\o-r\ut\<\phi\o,r\uo\>)\tens(\phi\t-r'\ut\<\phi\t,r'\uo\>)
=(\theta\tens\theta)\delta\phi}
using the Lie cobracket $\delta_{\rm bos}$ on $\und\cg^*\lbiprod\cg$ from
Theorem~3.5. The braided-Lie cobracket of $\und\cg^*$ coincides with that of
$\cg^*$, i.e. $\und\delta\phi=\delta\phi$. We also use the quasitriangular form
of $\cg$ to compute its Lie cobracket on $r\ut$.

Note that another way to present the  result is that $\pi(\xi)=\xi$ and
$\pi(\phi)=-r\ut\<\phi,r\uo\>$ is a Lie bialgebra projection $D(\cg)\to\cg$
split by  the inclusion of $\cg$, and recognise $\und\cg^*$ as the image under
$\theta$ of the braided-Lie bialgebra kernel of this according to Theorem~3.7.
The computations involved are similar to the above proofs for $\theta$. Similar
formulae are obtained if one takes $\pi(\phi)=r\uo\<\phi,r\ut\>$, corresponding
to transmutation with respect to the conjugate quasitriangular structure.
\endproof

This is the Lie version of the result for the quantum double of a
quasitriangular Hopf algebra in \cite{Ma:dou}. It completes the partial result
in \cite{Ma:skl} where, in the absence of a theory of braided-Lie bialgebras
we could only give the result $D(\cg)\isom \und\cg\lcross\cg$ in the
factorisable case (where $\und\cg^*\isom\cg$) and only as a Lie algebra
isomorphism. Since $\und\cg\lcross\cg$ by $\ad$ is easily seen to be isomorphic
to a direct sum Lie algebra $\cg\oplus\cg$,  one recovers the result that
$D(\cg)$ in the factorisable case is a Lie algebra direct sum, but now with a
certain Lie bialgebra structure (namely the double cross cosum
$\cg\codcross\cg$ in \cite{Ma:book}).

More recently, we have obtained a more general `double bosonisation'
theorem\cite{Ma:dbos} which yields as output quasitriangular Hopf algebras. It
provides an inductive construction for factorisable quasitriangular Hopf
algebras such as $U_q(\cg)$. The Lie version of this is as follows.  We suppose
$\cc,\cb$ are dually paired in the sense of a morphism $\<\ ,\
\>:\cc\tens\cb\to k$ such that the Lie bracket of one is adjoint to the Lie
cobracket of the other, and vice versa. The nicest case is where $\cb$ is
finite-dimensional and $\cc=\cb^*$ as in Lemma~3.4, but we do not need to
assume this for the main construction.

\begin{theorem} For dually paired braided Lie bialgebras
$\cb,\cc\in{}_{\cg}\CM$ the vector space $\cb\oplus\cg\oplus\cc$ has a unique
Lie bialgebra structure $\cb\lbiprod\cg\rbiprod\cc^{\rm
op}$, the {\em double-bosonisation}, such that $\cg$ is a sub-Lie bialgebra,
$\cb,\cc^{\rm op}$ are sub-Lie algebras, (\ref{liebos}) and
\cmath{{}[\xi,x]=\xi\la x,\quad [\xi,\phi]=\xi\la \phi,\quad
[x,\phi]=x\Bo\<\phi,x\Bt\>+\phi\Bo\<\phi\Bt,\phi\>+2r_+\uo\<\phi,r_+\ut\la
x\>\\
\delta x=\und\delta x+r\ut\tens r\uo\la x-r\uo\la x\tens r\ut,\quad
\delta\phi=\und\delta\phi+r\ut\la \phi\tens r\uo-r\uo\tens r\ut\la \phi}
$\forall x\in\cb,\xi\in\cg$ and $\phi\in\cc$. Here $\und\delta
x=x\Bo\tens x\Bt$.
\end{theorem}
\proof  Here $\cb,\cg$ clearly form the bosonisation Lie bialgebra
$\cb\lbiprod\cg$ from Theorem~3.5. In the same way, we recognise  $\cg\rbiprod
\cc^{\rm op}$ as the bosonisation of $\cb^{*\rm op}$ as a braided-Lie bialgebra
in the category of $\cg$-modules with opposite infinitesimal braiding (see the
remark below Lemma~3.4). Since these are already known to form Lie bialgebras,
the coJacobi identity for the double-bosonisation holds, as well as the
1-cocycle axiom for all cases except $\delta([x,\phi])$ mixing $\cb,\cc$. We
outline the proof of this remaining case. From the definition of
$\cb\lbiprod\cg\rbiprod\cc^{\rm op}$, we have
\align{\delta([x,\phi])\equad
&&=\delta(x\Bo\<\phi,x\Bt\>+\phi\Bo\<\phi\Bt,x\>+2r_+\uo\<\phi,r_+\ut\la x\>)\\
&&=x\Bo\Bo\tens x\Bo\t\<\phi,x\Bt>+r\ut\tens r\uo\la x\Bo\<\phi,x\Bt\>-r\uo\la
x\Bo\tens r\ut\<\phi,x\Bt\>\\
&&\quad +\phi\Bo\Bo\tens\phi\Bo\Bt\<\phi\Bt,x\>+r\ut\la\phi\Bo\tens
r\uo\<\phi\Bt,x\>-r\uo\tens r\ut\la\phi\Bo\<\phi\Bt,x\>\\
&&\quad+2\delta r_+\uo\<\phi,r_+\ut\la x\>\\
\ad_x\delta\phi\equad &&=[x,\phi\Bo]\tens
\phi\Bt+\phi\Bo\tens[x,\phi\bt]+[x,r\ut\la \phi]\tens r\uo\\
&&\quad +r\ut\la\phi\tens [x,r\uo]-[x,r\uo]\tens
r\ut\la\phi-r\uo\tens[x,r\ut\la\phi]\\
&&=x\Bo\<\phi\Bo,x\Bt\>\tens\phi\Bt+\phi\Bo\Bo\<\phi\Bo\Bt,x\>\tens\phi\Bt
+2r_+\uo\<\phi\Bo,r_+\ut\la x\>\tens\phi\Bt\\
&&\quad+\phi\Bo\tens
x\Bo\<\phi\Bt,x\Bt\>+\phi\Bo\tens\phi\Bt\Bo\<\phi\Bt\Bt,x\>+\phi\Bo\tens
2r_+\uo\<\phi\Bt,r_+\ut\la x\>\\
&&\quad -r\ut\la\phi\tens r\uo\la x+r\uo\la x\tens r\ut\la\phi+[x,r\ut\la
\phi]\tens r\uo-[x,r\uo]\tens r\ut\la\phi-r\uo\tens[x,r\ut\la\phi].}
In a similar way, one has
\align{-\ad_\phi\delta x\equad &&=x\Bo\Bo\<\phi,x\Bo\Bt\>\tens
x\Bt+\phi\Bo\<\phi\Bt,x\Bo\>\tens x\Bt+2r_+\uo\<\phi,r_+\ut\la x\t\>\\
&&\quad+x\Bo\tens
x\Bt\Bo\<\phi,x\Bt\Bt\>+x\Bo\<\phi\Bt,x\Bt\>\tens\phi\Bo+x\Bo\tens
2r_+\uo\<\phi,r_+\ut\la x\Bt\>\\
&&\quad+r\ut\la\phi\tens r\uo\la x-r\uo\la x\tens r\ut\la\phi+r\ut\tens[r\uo\la
x,\phi]-[r\uo\la x,\phi]\tens r\ut.}
Adding the latter two expressions and comparing with $\delta([x,\phi])$ we see
that the terms of the form $r\ut\la\phi\tens r\uo\la x$ etc. immediately
cancel, the terms of the form $\phi\Bo\<\phi\Bt,x\Bt\>\tens x\Bt$ etc.
(involving Lie cobrackets of both $x$ and $\phi$) cancel by antisymmetry of the
Lie cobrackets, and the terms of the form $x\Bo\Bo\<\phi,x\Bo\Bt\>\tens x\Bt$
etc. (involving iterated Lie cobrackets of either $x$ or $\phi$) cancel using
antisymmetry of the Lie cobrackets and the coJacobi identity
$(\id\tens\delta)\delta +{\rm cyclic}=0$ for $\cb$ and $\cc$. Hence the
1-cocycle identity for this case reduces to the more manageable
\align{&&\equad r\ut\tens r\uo\la x\Bo\<\phi,x\Bt\>+r\ut\la\phi\Bo\tens
r\uo\<\phi\Bt,x\>+[2r_+\uo,r'\uo]\tens r'\ut \<\phi,r_+\ut\>-{\rm flip}\\
&&=2r_+\uo\<\phi\Bo,r_+\ut\la x\>\tens\phi\Bt+2r_+\uo\<\phi,r_+\ut\la
x\Bo\>\tens x\Bt\\
&&\quad +[x,r\ut\la \phi]\tens r\uo-[r\uo\la x,\phi]\tens r\ut -{\rm flip}}
where `-flip' means to subtract all the same expressions with the opposite
tensor product. We used antisymmetry of the Lie cobrackets and the
quasitriangular of $\cg$ for $\delta r_+\uo$.
On then has to put in the stated definitions of the Lie brackets $[x,r\ut\la
\phi]$ and $[r\uo\la x,\phi]$ and use $\cg$-covariance of the pairing, and of
the braided-Lie brackets and cobrackets to obtain equality.

Note that by comparing the Lie bosonisation formulae with the braided group
case, we can read off the Lie  double-bosonisation formulae from the braided
group case given in the required left-module form in the appendix of
\cite{Ma:conf}. The only subtlety is that in the Lie case we can eliminate the
categorical pairing $\ev$ (corresponding to the categorical dual $\cb^\star$ in
the finite-dimensional case): $\cc,\cb$ are categorically paired by
$\ev:\cc\tens\cb\to k$ {\em iff} $\<\ ,\ \>=-\ev$ is a ($\cg$-equivariant)
ordinary duality pairing. Then one obtains the $[x,\phi]$ relations as stated.
Finally, in \cite{Ma:dbos} it is explicitly shown that the double-bosonisation
is built on the tensor product vector space. The analogous arguments now prove
that the Lie double bosonisation is built on the direct sum vector space.
\endproof

\begin{propos}
Let $\cb\in{}_{\cg}\CM$ be a finite-dimensional braided-Lie bialgebra with dual
$\cb^*$. Then the double-bosonisation $\cb\lbiprod\cg\rbiprod\cb^{*\rm op}$ is
quasitriangular, with
\[ r^{\rm new}=r+\sum_a f^a\tens e_a,\]
where $\{e_a\}$ is a basis of $\cb$ and $\{f^a\}$ is a dual basis, and $r$ is
the quasitriangular structure of $\cg$. If $\cg$ is factorisable then so is the
double-bosonisation.
\end{propos}
\proof We show first that the Lie cobracket of the double-bosonisation has the
form $\delta=\extd r_{\rm new}$. With summation over $a$ understood, we have
\align{&&\equad [\phi,r\uo_{\rm new}]\tens r\uo_{\rm new}+r\uo_{\rm
new}\tens [\phi,r\ut_{\rm new}]\\
&&=[\phi,r\uo]\tens r\ut+r\uo\tens[\phi,r\ut]-[\phi,f^a]_{\cb^*}\tens
e_a+f^a\tens[\phi,e_a]\\
&&=-[\phi,f^a]_{\cb^*}\tens e_a - r\uo\la\phi\tens r\ut - r\uo\tens
r\ut\la\phi\\
&&\quad -f^a\tens e_a\Bo\<\phi,e_a\Bt\>-f^a\tens\phi\Bo\<\phi\Bt,e_a\>-f^a\tens
2r_+\uo\<\phi,r_+\ut\la e_a\>\\
&&=\und\delta\phi- r\uo\la\phi\tens r\ut - r\uo\tens
r\ut\la\phi+2r_+\ut\la\phi\tens r_+\uo=\delta\phi}
as required. Here $\<[f^a,\phi]_{\cb^*},x\>\tens e_a=f^a\tens
e_a\Bo\<\phi,e_a\Bt\>$ since both evaluate against $x\in\cb$ to
$x\Bo\<\phi,x\Bt\>$. The suffix $\cb^*$ is to avoid confusion with the Lie
bracket inside the double-bosonisation, which is that of $\cb^{*\rm op}$ on
these elements.
Similarly,
\align{&&\equad [x,r\uo_{\rm new}]\tens r\uo_{\rm new}+r\uo_{\rm new}\tens
[x,r\ut_{\rm new}]\\
&&=-r\uo\la x+r\ut-r\uo\tens r\ut\la x+[x,f^a]\tens e_a+f^a\tens[x,e_a]\\
&&=-r\uo\la x+r\ut-r\uo\tens r\ut\la x+f^a\tens[x,e_a]+
x\Bo\<f^a,x\Bt\>\tens e_a+f^a\Bo\<f^a\Bt,x\>\tens e_a\\
&&\quad+2r_+\uo\tens r_+\ut\la x\\
&&=\und\delta x - r\uo\la x\tens r\ut+r\ut\tens r\uo\la x=\delta x.}
Here $f^a\Bo\<f^a\Bt,x\>\tens e_a=-f^a\tens[x,e_a]$ as both evaluate against
$\phi\in\cb^*$ to $\phi\Bo\<\phi\Bt,x\>$. Since the Lie cobracket of the
double-bosonisation is antisymmetric, we
conclude also that $2r_{+\rm new}$ is $\ad$-invariant.

Finally, we verify the CYBE for $r_{\rm new}$. Actually, once $\delta=\extd
r_{\rm new}$ has been established, the CYBE is equivalent to
$(\delta\tens\id)r_{\rm new}=r\uo_{\rm new}\tens r'\uo_{\rm new}\tens
[r\ut_{\rm new},r'\ut_{\rm new}]$ (see \cite{Ma:book}). Note that
\[ \und\delta f^a\tens e_a=f^a\tens f^b\tens [e_a,e_b]\]
(sum over $a,b$) since  evaluation against $x,y\in\cb$ gives $[x,y]$ in both
cases. Then
\align{(\delta\tens\id)r_{\rm new}\equad&&=(\delta\tens\id)r+\delta f^a\tens
e_a\\
&&=r\uo\tens r'\uo\tens [r\ut,r'\ut]+f^a\tens f^b\tens [e_a,e_b]+r\ut\la
f^a\tens r\uo\tens e_a-r\uo\tens r\ut\la f^a\tens e_a\\
&&= r\uo\tens r'\uo\tens [r\ut,r'\ut]+f^a\tens f^b\tens [e_a,e_b]-f^a\tens
r\uo\tens [r\ut,e_a]+r\uo\tens   f^a\tens [r\ut, e_a]}
as required. We used $\cg$-covariance of the pairing, so that $\xi\la f^a\tens
e_a=-f^a\tens\xi\la e_a=-f^a\tens[\xi,e_a]$ for all $\xi\in\cg$.

If $\cg$ is factorisable then $2r_{+{\rm new}}$ as a map
$(\cb\lbiprod\cg\rbiprod\cb^{*\rm op})^*\to \cb\lbiprod\cg\rbiprod\cb^{*\rm
op}$ has $\cg$ in its image, by restricting to $\cg$. It has $\cb$ in its image
by restricting to $\cb$, and $\cb^*$ in its image by restricting to $\cb^*$. So
the double-bosonisation is again factorisable. \endproof

There is also a more general double-bisum construction
$\cb\lbiprod\cf\,\rbiprod\cc^{\rm op}$ containing
biproducts $\cb\lbiprod\cf$ and $\cf\,\rbiprod\cc^{\rm op}$ (with
$\cc,\cb\in{}^{\cf}_{\cf}\CM$ suitably paired braided-Lie bialgebras) and
reducing to the double-bosonisation in the case when $\cc,\cb$ are in the image
of the functor in Lemma~3.8.

Double bosonisation reduces to Drinfeld's double $D(\cb)$ when $\cf=0$ (then a
braided-Lie bialgebra reduces to an ordinary Lie bialgebra). And because it
preserves factorisability, it provides an inductive construction for new
factorisable quasitriangular Lie bialgebras from old ones. We will see in the
next section that it can be used as a coordinate free version of the idea of
adjoining a node to a Dynkin diagram
(adjoining a simple root vector in the Cartan-Weyl basis). Moreover, building
up $\cg$ iteratively like this also builds up the quasitriangular
structure~$r$.

\section{Parabolic Lie bialgebras and Lie induction}

In this section we give some concrete examples and applications of
the above theory. We work over $\C$. We begin with the simplest
example of a braided-Lie bialgebra, with zero Lie bracket and zero
Lie  cobracket. According to Definition~2.2 this means precisely
modules of our background quasitriangular Lie bialgebras for which
the infinitesimal braiding cocycle $\psi$ vanishes.

\begin{propos} Let $\cg$ be a semisimple factorisable (s.s.f) Lie
bialgebra
 and $\cb$ an isotypical representation such that $\Lambda^2\cb$ is
isotypical. Then $\cb$ with zero bracket and zero cobracket is a
braided-Lie
bialgebra in ${}_{\tilde{\cg}}\CM$, where $\tilde{\cg}$ is a central
extension.
\end{propos}
\proof Let $c=r_+\uo r_+\ut$ in $U(\cg)$. Since $r_+$ is
ad-invariant, $c$ is central. Moreover, $2r_+=\Delta c-(c\tens
1+1\tens c)$ where $\Delta$ is the coproduct of $U(\cg)$ as a Hopf
algebra.
Since $\cb$ is assumed isotypical, the action of $c$ on it is by
multiplication by a scalar, say $\lambda_1$. Since $\Lambda^2\cb$ is
assumed isotypical, the action of $c$ on it, which is the action of
$\Delta c$ in each factor, is also multiplication by a scalar, say
$\lambda_2$. Then $\psi(x\tens y)=(\Delta c-(c\tens 1+1\tens
c))\la(x\tens y-y\tens x)=(\lambda_2-2\lambda_1)(x\tens y-y\tens
x)=\lambda(x\tens y-y\tens x)$ say, where $\lambda$ is a constant.

Now, $\cb$ with the zero bracket and cobracket is not a braided group
in ${}_{\cg}\CM$ unless our cocycle $\psi$ vanishes. However, in the
present case we can neutralise the cocycle with a central extension.
Thus, let $\tilde{\cg}=\C\oplus\cg$ with $\C$ spanned by $\dila$,
say.
We take the Lie bracket, quasitriangular structure and Lie cobracket
\[ [\xi,\dila]=0,\quad \tilde r=r-{\lambda\over
2}\dila\tens\dila,\quad \delta\dila=0\]
for all $\xi\in\cg$. In this way, $\tilde{\cg}$ becomes a
quasitriangular Lie bialgebra.
We consider $\cb\in{}_{\cg}\CM$ by $\dila\la x=x$ for all $x\in\cb$.
The infinitesimal braiding on $\cb$ in this category is
$\tilde\psi(x\tens y)=2\tilde r_+\la(x\tens y-y\tens x)=\psi(x\tens
y)-\lambda(x\tens y-y\tens x)=0$. So $\cb$ is a braided-Lie bialgebra
in this category. \endproof

The constant $\lambda$ is the infinitesimal analogue of the so-called
quantum group normalisation constant. The central extension is the
analogue of the central extension by a `dilaton' needed for the
quantum planes to be viewed as braided groups\cite{Ma:poi}. We see
now the infinitesimal
analogue of this phenomenon.

Next, we can apply Theorem~3.5 and
obtain a Lie bialgebra $\cb\lbiprod\tilde{\cg}$ as the bosonisation of
$\cb$. Moreover, double-bosonisation provides a still bigger and
factorisable Lie algebra containing $\cb\lbiprod\tilde{\cg}$.

\begin{corol} Let $\cg$ be simple and strictly quasitriangular, and
$\cb$ a finite-dimensional irreducible representation with
$\Lambda^2\cb$ isotypical. Then
the double bosonisation $\cb\lbiprod\tilde{\cg}\rbiprod\cb^*$ from
Theorem~3.10 is again simple, strictly quasitriangular and of
strictly greater rank.
\end{corol}
\proof The Lie bracket in the double-bosonisation in Theorem~3.10 and
the form of $\tilde r$ are
\[ {}[\xi,x]=\xi\la x,\quad [\dila,x]=x,\quad
[\xi,\phi]=\xi\la\phi,\quad [\dila,\phi]=-\phi,\quad [\xi,\dila]=0\]
\[{}[x,\phi]=2r_+\uo\<\phi,r_+\ut\la x\>-\lambda\dila\<\phi,x\>\]
for all $\xi\in\cg$, $x\in\cb$ and $\phi\in\cg^*$.
Consider $I\subseteq \cb\oplus\cg\oplus\C\oplus\cb^*$ an ideal of the
double-bosonisation. Let $I_{\cb},I_{\cb^*},I_{\cg},I_{\C}$ be the
components of $I$ in the direct sum. By the relation $[\xi,x]=\xi\la
x$, $I_{\cb}\subseteq\cb$ is a subrepresentation under $\cg$. Since
$\cb$ is irreducible, $I_{\cb}$ is either zero or $\cb$. Similarly
for $I_{\cb^*}$. Likewise $I_{\cg}$ is zero or $\cg$ as $\cg$ is
simple. Finally, $I_{\C}$ is zero or $\C$ by linearity. We therefore
have 16 possibilities to consider for whether $\C,\cg,\cb,\cb^*$ are
contained or not in $I$.  (i) if $\cg$ is contained, then since $\cb$
is irreducible, the relation $[\xi,x]=\xi\la x$ spans $\cb$ for any
fixed $x$, and hence is certainly not always zero. So $\cb$ is
contained, and likewise $\cb^*$ is contained if $\cg$ is. In this
case, the $[x,\phi]$ relation means that $\C$ is contained and $I$ is
the whole space. (ii) if  $\cb$ is contained then the $[x,\phi]$
relation and $2r_+$ nondegenerate means that $\cg$ and $\C$ are
contained and hence $I$ is the whole space. (iii) Similarly if
$\cb^*$ is contained. (iv) Finally, if $\C$ is contained then the
relation $[\dila,x]=x$ implies that $\cb$ is contained and
hence $I$ is the whole space. Hence $I$ is zero or the whole space,
as required. The new quasitriangular structure is non-zero since its
component in $\cg\tens\cg$ is non-zero. The rank is clearly increased
by at least 1 due to the addition of $\dila$. \endproof

Thus the double-bosonisation in Theorem~3.10 provides an inductive
construction for simple strictly quasitriangular Lie bialgebras.
It is possible to see that the fundamental
representations of $su_n$ or
$so_n$ take us up to $su_{n+1}$ and $so_{n+1}$, i.e. precisely take
us up the ABD series
in the usual classification of Lie algebras. Moreover, we see the
role of the single bosonisation in Theorem~3.5:

\begin{example} Consider $\cg=su_2$ with the Drinfeld-Sklyanin
quasitriangular structure. The 2-dimensional irreducible
representation $\cb$ is a braided-Lie bialgebra via Proposition~4.1.
Its bosonisation
$\C^2\lbiprod\widetilde{su_2}$ is the maximal parabolic of the
double bosonisation
$\C^2\lbiprod\tilde{\cg}\rbiprod\C^2=su_3$. Explicitly,
it is the Lie algebra of $su_2$ and
\cmath{{} [x,y]=0,\quad [X_+,x]=0,\quad [X_+,y]=x,\quad
[X_-,x]=y,\quad [X_-,y]=0\\
{}[H,x]=x,\quad [H,y]=-y,\quad [\dila,H]=0,\quad
[\dila,X_\pm]=0,\quad [\dila,x]=x,\quad [\dila,y]=y}
where $\{x,y\}$ are a basis of $\C^2$ and $H,X_\pm$ are the standard
$su_2$ Chevalley generators. The Lie cobracket on the generators is
\[\delta\dila=0,\quad\delta X_\pm=\h X_\pm\wedge H,\quad \delta
x=\h x\wedge
h\]
where $h=-{1\over 2}H-{3\over 2}\dila$ and $
\wedge=(\id-\tau)\circ\tens$.
\end{example}
\proof Note that we work over $\C$, but there is a natural real
forms justifying the notation.  Here $\Lambda^2\cb$ is the
1-dimensional (i.e. spin 0) representation of $su_2$. The standard
quasitriangular structure of $su_2$ is
\[ r={1\over 4}H\tens H+X_+\tens X_- \]
Then $c=r_+\uo r_+\ut$ is twice the quadratic Casimir in its usual
normalisation. Hence its value in the $(2j+1)$ dimensional (i.e. spin
$j$) irreducible representation is $j(j+1)$. In the present case, we
have $\lambda=0.(0+1)-2.\h(\h+1)=-{3\over 2}$ in Proposition~4.1. We
therefore make the central extension to $\tilde{\cg}$ and apply
Theorem~3.5. The Lie algebra of the bosonisation is the action of
$\widetilde{\cg}$. Its explicit form in the representation
$\rho(X_+)=\pmatrix{0&1\cr0&0}$, $\rho(X_-)=\pmatrix{0&0\cr 1&0}$ and
$\rho(H)=\pmatrix{1&0\cr0&-1}$ is $X_+\la x=0$, $X_+\la y=x$, $X_-\la
x=y$, $X_-\la y=0$, $H\la x=x$ and $H\la y=-y$, giving the Lie
bracket stated. The Lie cobracket is $\delta x=0+\tilde
r\ut\wedge\tilde r\uo\la x={1\over 4}H\wedge H\la x+{3\over
4}\dila\wedge\dila\la x=\h x\wedge h$ as stated. We identify
$X_\pm=X_{\pm1},H=H_1$ as a sub-Lie
algebra of $su_3$ and $x=X_{-2},h=H_2$ as the remaining Chevalley
generators of its standard maximal parabolic. Finally, let $\cb^*$
have dual basis $\{\phi,\psi\}$.
By a similar computation to the above, we obtain
$\tilde{su_2}\rbiprod \C^2$ with Lie bracket
\cmath{{}[\phi,\psi]=0,\quad [X_+,\phi]=-\phi,\quad
[X_+,\psi]=0,\quad [X_-,\phi]=0,\quad [X_-,\psi]=-\phi\\
{}[H,\phi]=-\phi, \quad [H,\psi]=\psi,\quad
[\dila,\phi]=-\phi,\quad[\dila,\psi]=-\psi.}
Among the further $\cb,\cb^*$ brackets in the double bosonisation in
Theorem~3.10, we have
$[x,\phi]=2r_+\uo\<\phi,r_+\ut\la x\>+{3\over 2}\dila\<\phi,x\>=\h
H\<\phi,H\la x\>+0+{3\over 2}\dila\<\phi,x\>=-h$.
{}From these relations we find that $\phi=X_{+2}$ and $\psi=X_{+12}$
explicitly identifies the double bosonisation as $su_3$. The Lie
cobracket on $\phi$ is $\delta\phi=\tilde r\ut\la\phi\wedge \tilde
r\uo={1\over 4}H\la\phi\wedge H+{3\over
4}\dila\la\phi\wedge\dila=\h\phi\wedge h$. This conforms with the
standard Lie cobracket for $su_3$. Indeed, the quasitriangular
structure of the double bosonisation in Theorem~3.10   reproduces the
Drinfeld-Sklyanin quasitriangular structure of $su_3$. \endproof

This is far from the only braided-Lie bialgebra in the category of
$\widetilde{su_2}$-modules, however.

\begin{example} Consider $\cg=su_2$ with the Drinfeld-Sklyanin
quasitriangular structure. The 3-dimensional irreducible
representation $\cb$ is a braided-Lie bialgebra via Proposition~4.1.
Its bosonisation $\R^3\lbiprod\widetilde{so_3}$ is the maximal
parabolic of the double bosonisation $so_5$. Explicitly, it is the
Lie
algebra of $so_3$ and
\cmath{{}[x_i,x_j]=0,\quad [e_i,x_j]=\sum_k\eps_{ijk}x_k,\quad
[\dila,x_j]=x_j }
where $i,j,k=1,2,3$ and $\eps$ is the totally antisymmetric tensor
with $\eps_{123}=1$. Here $\{e_i\}$ are the vector basis of $so_3$.
The Lie
cobracket is
\cmath{\delta e_1=\imath e_1\wedge e_3,\quad \delta e_2=\imath
e_2\wedge e_3,\quad\delta e_3=0,\quad\delta\dila=0,\quad\delta
x_1=(\imath
e_1+e_2)\wedge x_3 +x_2\wedge e_2 +\dila\wedge x_1\\
\delta x_2=x_3\wedge(e_1-\imath e_2) +e_3\wedge x_1+\dila\wedge
x_2,\quad \delta x_3=(e_1-\imath e_2)\wedge x_2+x_1\wedge (\imath
e_1+e_2)+\dila\wedge x_3}
\end{example}
\proof Here $\Lambda^2\cb$ is also the 3-dimensional (i.e. spin 1)
representation. Hence, from the first part of the proof of
Example~4.2, we have $\lambda=1.(1+1)-2.1.(1+1)=-2$. The Lie algebra
$so_3$ in the vector basis is $[e_1,e_2]=e_3$ and
cyclic rotations of this, and the Drinfeld-Sklyanin quasitriangular
structure in this basis is \cite[Ex. 8.1.13]{Ma:book}
\[ r=-e_i\tens e_i+\imath(e_1\tens e_2-e_2\tens e_1)\]
We add $\dila\tens\dila$ to give the quasitriangular structure
$\tilde r$. The action on $\C^3$ with basis $x_i$ is $[e_1,x_2]=x_3$
and cyclic rotations of this. This immediately provides the Lie
algebra of the bosonisation. The Lie cobracket from Theorem~3.5 is
\[ \delta x_i=\imath e_2\wedge [e_1,x_i]-\imath e_1\wedge
[e_2,x_i]+e_j\wedge \eps_{ijk}x_k+\dila\wedge x_i\]
with computes as stated. \endproof

This example is manifestly the Lie  algebra of motions plus dilation
of $\R^3$ as a sub-Lie algebra of the conformal Lie algebra
$so(1,4)$, equipped now with a Lie bialgebra structure. At the level
of complex Lie algebra, it is the maximal parabolic of $so_5$.
The generator $\dila$ is called the `dilaton' in the corresponding
quantum groups literature. We likewise obtain natural maximal
parabolics for the whole ABD series by bosonisation of the
fundamental representation $\cb$.

On the other hand, these steps for
other Lie algebras
can involve less trivial braided-Lie bialgebras $\cb$ (with non-zero
bracket
and cobracket). The general case is as follows. We consider simple Lie
algebras $\cg$ associated to root systems in the usual conventions.
Positive roots are denoted $\alpha$, with length $d_\alpha$. The
Cartan-Weyl basis has root vectors $X_{\pm\alpha}$ and Cartan
generators $H_i$ corresponding to the simple roots $\alpha_i$. We
define $d_\alpha H_\alpha  =\sum_i n_i d_i H_i$ if $\alpha=\sum_i
n_i\alpha_i$. We take the Drinfeld-Sklyanin quasitriangular structure
in its general form
\eqn{driskl}{ r=\sum_\alpha d_\alpha X_\alpha\tens
X_{-\alpha}+\h\sum_{ij} A_{ij}H_i\tens H_j,}
where $A_{ij}=d_i(a^{-1})_{ij}$. Here $a$ is the Cartan matrix.
The corresponding Lie cobracket is $\delta X_{\pm i}=\h X_{\pm
i}\wedge H_i$ and $\delta H_i=0$ on the generators.

\begin{propos} Let $i_0$ be a choice of simple root such that its
deletion again generates the root system of a simple Lie algebra,
$\cg_0$. Let $\cb_-\subset\cg$ be the standard (negative) Borel and
let $\cf\subset\cb_-$ denote the sub-Lie algebra excluding all
vectors generated by $X_{-i_0}$. Both $\cb_-$ and $\cf$ are sub-Lie
bialgebras of $\cg$ and
\[ \cb_-{{\buildrel \pi\over \to}\atop\hookleftarrow} \cf,\quad
\pi(H_i)=H_i,\quad
\pi(X_{-\alpha})=\cases{0&if $\alpha$ contains $\alpha_{i_0}$\cr
X_{-\alpha}&else}\]
is a split Lie bialgebra projection.  Then
$\cb=\ker\pi$ is the Lie ideal generated by $X_{-i_0}$ in $\cb_0$ and
is a braided-Lie bialgebra in
${}_{\cf}^{\cf}\CM$ by Theorem~3.7.
\end{propos}
\proof Here $\cf$ is generated by all the $H_i$ and  only those $X_{-j}$ where
$j\ne i_0$, i.e. spanned by the $H_i$ and $\{X_{-\alpha}\}$ such that
$\alpha$ does not contain $\alpha_{i_0}$. It is clearly a sub-Lie
algebra of $\cb_-$. We show first that it is a sub-Lie bialgebra.
First of all, note that the Lie coproduct in $\cg$ has the general
form
\[ \delta X_{\pm\alpha}={d_\alpha\over 2} X_{\pm\alpha}\wedge H_\alpha +
\sum_{\beta+\gamma=\alpha}c_{\pm\beta,\pm\gamma}X_{\pm\beta}\wedge X_{\pm
\gamma}\]
where the sum is over positive root $\beta,\gamma$ adding up to
$\alpha$ and the $c$ are constants. The proof is by induction (being careful
about signs). From
the Lie bialgebra cocycle axiom and the induction hypothesis,
\align{\delta([X_i,X_\alpha])\equad&&={d_\alpha\over
2}[X_i,X_\alpha]\wedge H_\alpha+{d_\alpha\over 2}X_\alpha\wedge
[X_i,H_\alpha]+\sum_{\beta+\gamma=\alpha}c_{\beta,\gamma}[X_i,X_\beta
]\wedge X_\gamma\\
&&\quad+\sum_{\beta+\gamma=\alpha}c_{\beta,\gamma}X_\beta\wedge
[X_i,X_\gamma]-{d_i\over 2}[X_\alpha,X_i]\wedge H_i-{d_i\over
2}X_i\wedge[X_\alpha,H_i]\\
&&={d_{\alpha+\alpha_i}\over 2}[X_i,X_\alpha]\wedge
H_{\alpha+\alpha_i}+\alpha(d_iH_i)X_i\wedge
X_\alpha+\sum_{\beta+\gamma=\alpha}c_{\beta,\gamma}[X_i,X_\beta]
\wedge X_\gamma\\
&&\quad+\sum_{\beta+\gamma=\alpha}c_{\beta,\gamma}X_\beta\wedge
[X_i,X_\gamma]}
if $\alpha+\alpha_i$ is a positive root. We used the identities
$[d_\alpha H_\alpha,X_i]=\alpha(d_iH_i)X_i$ and
$[d_iH_i,X_\alpha]=\alpha(d_iH_i)X_\alpha$. Since all positive root
vectors are obtained by iterated Lie brackets of the $X_i$, we
conclude the result (the argument for negative roots is similar).

{}From this form, it is clear first of all that $\delta$ restricts to
$\cb_-\to\cb_-\tens\cb_-$, so this becomes a sub-Lie bialgebra of
$\cg$ (this is well-known). Moreover, if $\alpha$ does not involve
$\alpha_{i_0}$ then neither can positive $\beta,\gamma$ such that
$\beta+\gamma=\alpha$. Hence $\cf$ is a Lie sub-bialgebra of $\cb_-$.

Finally, $\pi$ is clearly a Lie algebra map by considering the cases
separately. For elements of $\cf\tens\cf$ we know that $\pi\circ[\ ,\
]=[\ ,\ ]=[\pi(\ ),\pi(\ )]$ since $\cf$ is closed, while if $\alpha$
involves $\alpha_{i_0}$ then so does $\alpha+\beta$  and
$\pi([X_{-\alpha},X_{-\beta}])=0=[\pi(X_{-\alpha}),\pi(X_{-\beta})]$.
Moreover, $\pi$ is a Lie coalgebra map on $\cf$ since
$\delta\cf\subset\cf\tens\cf$ as shown above. Finally,
$(\pi\tens\pi)\delta X_{-i_0}=0=\delta\pi (X_{-i_0})$ from the simple
form of $\delta$ on the generators.

Therefore we may apply Theorem~3.7 and obtain a braided-Lie bialgebra
$\cb=\ker\pi$. Here $\cb\subset\cb_-$ is the Lie ideal generated
by $X_{-{i_0}}$, i.e. spanned by $\{X_{-\alpha}\}$ where $\alpha$
contains $\alpha_{i_0}$.

The braided-Lie cobracket of $\cb$ from Theorem~3.7 is $\und\delta
X_{-\alpha}=\sum c_{-\beta,-\gamma}X_{-\beta}\wedge X_{-\gamma}$, the
part of the Lie cobracket $\delta X_{-\alpha}$ in which both
$\beta,\gamma$ contain $\alpha_{i_0}$. The action of $\cf$ is by Lie
bracket in $\cb_-$ and the Lie coaction of $\cf$ is
$\beta(X_{-\alpha})=-{d_\alpha\over 2}H_\alpha\tens X_{-\alpha}+\sum
c_{-\beta,-\gamma}X_{-\beta}\tens X_{-\gamma}$ where the sum is the
part of $\delta X_{-\alpha}$ where $\beta$ does not contain
$\alpha_{i_0}$. \endproof

This constructs the required braided-Lie bialgebra for the general
case. Although obtained here in the category of $D(\cf)$-modules, this action
is
compatible with an action of the central extension $\tilde{\cg_0}\subset\cg$.
It is easy to see that there is a unique element $\dila\in\cg$,
$\dila\notin\cg_0$ which commutes with the image of $\cg_0$. It is determined
by the Cartan matrices
of $\cg,\cg_0$. Viewed in $\cg$, this $\tilde{\cg_0}$ acts on $\cg$ by the
adjoint action and
this action restricts to $\cb$. In this way, $\cb$ becomes a braided-Lie
bialgebra in
${}_{\tilde{\cg_0}}\CM$.

\begin{example} When $\cg=g_2$ and $\cg_0=su_2$, we obtain the 5-dimensional
braided-Lie bialgebra where $su_2$ acts  as the $4\oplus 1$ dimensional (i.e.
the spin $3\over 2$ and spin 0 representations). Both the Lie bracket and the
Lie cobracket are not identically zero.
\end{example}
\proof We take the Cartan matrix for $\cg$ as $\pmatrix{2&-1\cr -3&2}$. We
$i_0=1$ so that the required $su_2$ is spanned by $H_2,X_{\pm2}$. The negative
roots vectors $X_{-1},X_{-21},X_{-221},X_{-2221}$ span the 4-dimensional
representation of $su_2$, the  eigenvalues of the adjoint action of $\h H_2$
being $-{3\over 2},-{1\over 2},{1\over 2},{3\over 2}$ respectively. These and
the remaining negative root vector $X_{-12221}$ (which forms a $1$-dimensional
trivial representation of $su_2$) are a basis of
$\cb$. We then restrict the Lie bracket to $\cb$, the only non-zero entries
being
\[ [X_{-1},X_{-2221}]=X_{-12221}=[X_{-221},X_{-21}].\]
This is a central extension (by a cocycle) of the zero bracket on the
4-dimensional representation.
Lie cobracket can the be computed by projection of the Lie cobracket in $g_2$.
Since (as one may easily verify) the infinitesimal braiding is nontrivial, both
the braided-Lie bracket and braided-Lie cobracket on $\cb$ are not identically
zero. The element $\dila=2H_1+H_2$ commutes with $su_2$ and acts as the
identity in the 4-dimensional part of $\cb$.

\begin{example} When $\cg=sp_6$ and $\cg_0=sp_4$, we obtain the 5-dimensional
braided-Lie
bialgebra where $sp_4$ acts in the $4\oplus 1$ dimensional representation. Here
the 4-dimensional representation is the fundamental one of $sp_4$. Both the Lie
bracket and the Lie cobracket are not identically zero.
\end{example}
\proof We take the Cartan matrices for $\cg$ and $\cg_0$ as
\[ a=\pmatrix{2&-1&0\cr -1&2&-2\cr 0&-1&2},\quad a_0=\pmatrix{2&-2\cr -1&2}\]
where $i_0=1$. We identify $sp_4=C_2$ inside $sp_6$ as the root vectors
$X_{\pm2},X_{\pm3},X_{\pm23}$
and Cartan vectors $H_2,H_3$. The negative root vectors
$X_{-1},X_{-12},X_{-123},X_{-1223}$ form the 4-dimensional representation of
$sp_4$. These and the remaining negative root vector $X_{-11223}$  (which forms
a 1-dimensional trivial representation of $sp_4$) are a basis of $\cb$. We then
restrict the Lie bracket to $\cb$ and find that this is again a cocycle central
extension of the zero Lie bracket on the 4-dimensional representation, the only
non-zero entries being
\[ [X_{-12},X_{-123}]=\h X_{-11223},\quad [X_{-1},X_{-1223}]=X_{-11223}.\]
The infinitesimal braiding and the Lie cobracket are also nontrivial, as one
may verify by further computation. The element $\dila=H_1+H_2+H_3$ commutes
with $sp_4$ and acts as the identity in the 4-dimensional part of $\cb$.
\endproof

These examples show that the general case need not depart too far from the
setting of Proposition~4.1 and Corollary~4.2; we need to make a central
extension of the underlying irreducible representation to define $\cb$. By
construction, $\cb\lbiprod\tilde{\cg_0}$ is once again the maximal parabolic of
$\cg$ associated to $\alpha_{i_0}$. A similar construction works for more roots
missing, giving
non-maximal parabolics of the double-bosonisation We simply define $\pi$
setting to zero all the root vectors containing the roots to be deleted in
defining $\cg_0$.

\section{Concluding remarks}

We have given here the basic theory of braided-Lie algebras, obtained by
infinitesimalising the
existing theory of braided groups. We also outlined in Section~4 its
application to the inductive construction of
simple Lie algebras with their standard quasitriangular structures. Further
variations of these constructions are certainly possible, and by making them
one should be able to also obtain the other strictly quasitriangular Lie
bialgebras structures in the Belavin-Drinfeld classification\cite{BelDri:sol}.
For example, there is a twisting theory of quantum groups\cite{Dri:qua} and
braided groups\cite{Ma:qsta}. An infinitesimal version of the latter would
allow one to introduce additional twists at each stage of the inductive
construction of the simple Lie algebra.

Also, although we have (following common practice) named our Lie algebras by
their natural real forms, our Lie algebras in Section~4 were complex ones.
There is a theory of $*$-braided groups (real forms of braided groups) as well
as their corresponding bosonisations and
double-bosonisations\cite{Ma:qsta}\cite{Ma:conf}. The infinitesimal version of
these should yield, for example, $so(1,4)$ as a real form arising from the
double-bosonisation of the 3-dimensional  braided-Lie bialgebra in Example~4.4.
The construction of natural compact  real forms and the classification of real
forms would be a further goal. These are some directions for further work.

Finally, just as Lie bialgebras extend to Poisson-Lie groups, so
braided-Lie
bialgebra structures typically extend to the associated Lie group $B$ of
$\cb$, at least locally. First, one needs to exponentiate $\psi\in
Z^2_{\ad}(\cb,cb\tens\cb)$  to a
group cocycle $\Psi\in Z^2_{\Ad}(B,\cb\tens\cb)$. Since $\extd\delta=\psi$, we
should
likewise exponentiate $\delta$ to the group as a map $D:B\to\cb\tens\cb$ with
coboundary $\Psi$, and define from
this a `braided-Poisson bracket'. The latter will not, however, respect the
group
product in the usual way but rather up to a `braiding' obtained from $\psi$.
Details of these
braided-Poisson-Lie groups and  the example of the Kirillov-Kostant
braided-Poisson bracket from Example~3.3 extended to the group manifold (e.g.
to $SU_2$) will be presented elsewhere.


\end{document}